\documentclass[12pt]{article}
 \pdfoutput=1
\usepackage[cp1250]{inputenc}
\usepackage[OT1]{fontenc}

\usepackage{amsmath,amsfonts,amssymb,isomath,color}
\numberwithin{equation}{section}
\usepackage{mathrsfs} 
\usepackage[pdftex]{graphicx}

\usepackage{cite}
\usepackage{upgreek}
\usepackage[format=hang, indention=0cm, labelformat=simple, labelsep=space, textformat=simple, justification=justified, singlelinecheck=true, labelfont=bf, textfont=sl,width=1.0\textwidth,figurename=Figure]{caption}

\allowdisplaybreaks[1]

\newcommand{\dd}{\ensuremath{\text{d}}}
\newcommand{\cint}[1]{\mathchoice{-\hspace{-2.5ex}\int_{#1}}{-\hspace{-2.1ex}\int_{#1}}{{-}\hspace{-1.6ex}\int_{#1}}{-\hspace{-1.4ex}\int_{#1}}}
\newcommand{\fint}[1]{\mathchoice{=\hspace{-2.6ex}\int_{#1}}{=\hspace{-2.4ex}\int_{#1}}{{=}\hspace{-1.5ex}\int_{#1}}{=\hspace{-1.0ex}\int_{#1}}}

\newcommand{\Functional}[1]{\ensuremath{\mathscr #1}}
\newcommand{\Vector}[1]{\ensuremath{\mathbf #1}}
\newcommand{\VectorSymbol}[1]{\ensuremath{\mathbfit #1}}
\newcommand{\VectorMatrix}[1]{\ensuremath{\boldsymbol #1}}
\newcommand{\Matrix}[1]{\ensuremath{\boldsymbol {\mathsf #1}}}
\newcommand{\MatrixSymbol}[1]{\ensuremath{\matrixsym{#1}}}

\newtheorem{theorem}{Theorem}[section]

\newtheorem{remark}[theorem]{Remark}

\newcommand\calP{\mathscr{P}}
\newcommand\calPt{{\hspace{.1em}\widetilde{\hspace{-.1em}\mathscr{P}}}}
\renewcommand\AA{^{\mbox{\tiny A}}}
\newcommand\BB{^{\mbox{\tiny B}}}
\newcommand\AB{^{\mbox{\tiny AB}}}
\newcommand\ABast{^{\mbox{\tiny AB}\ast}}
\newcommand\C{_{\mbox{\tiny\rm C}}}
\newcommand\N{_{\mbox{\tiny\rm N}}}
\newcommand\D{_{\mbox{\tiny\rm D}}}
\newcommand\Dir{{\mbox{\tiny\rm D}}}

\newcommand\DN{_{\mbox{\tiny\rm DN}}}
\newcommand\CN{_{\mbox{\tiny\rm CN}}}
\newcommand\NN{_{\mbox{\tiny\rm NN}}}
\newcommand\DD{_{\mbox{\tiny\rm DD}}}
\newcommand\ND{_{\mbox{\tiny\rm ND}}}
\newcommand\NC{_{\mbox{\tiny\rm NC}}}
\newcommand\CD{_{\mbox{\tiny\rm CD}}}
\newcommand\DC{_{\mbox{\tiny\rm DC}}}
\newcommand\CC{_{\mbox{\tiny\rm CC}}}

\newcommand\GC{\Gamma\C}
\newcommand\GD{\Gamma\D}
\newcommand\GN{\Gamma\N}

\newcommand\OA{\Omega\AA}
\newcommand\OB{\Omega\BB}
\newcommand\GA{\Gamma\AA}
\newcommand\GB{\Gamma\BB}
\renewcommand\dot[1]{\mathchoice
                 {{\buildrel{\hspace*{.1em}\text{\LARGE.}}\over{#1}}}
                 {{\buildrel{\hspace*{.1em}\text{\Large.}}\over{#1}}}
                 {{\buildrel{\hspace*{.1em}\text{\large.}}\over{#1}}}
                 {{\buildrel{\hspace*{.1em}\text{\large.}}\over{#1}}}}

\newcommand\JUMP[2]{\mathchoice
                  {\big[\hspace*{-.3em}\big[#1\big]\hspace*{-.3em}\big]_{#2}^{}}
                   {[\hspace*{-.15em}[#1]\hspace*{-.15em}]_{#2}^{}}
                   {[\![#1]\!]_{#2}^{}}
                   {[\![#1]\!]_{#2}}^{}}
\newcommand{\wt}[1]{\mathchoice
     {\text{\small$\widetilde{\text{\normalsize$#1$}}\hspace*{.03em}$}}
                    {\text{\small$\widetilde{\text{\normalsize$#1$}}$}}
                    {\widetilde{#1\hspace*{-.02em}}\hspace*{.03em}}
                    {\widetilde{#1}}}
\renewcommand\Vector[1]{\ensuremath{#1}}
\renewcommand\VectorSymbol[1]{#1}
\newcommand{\weak}{\rightharpoonup}
\newcommand{\weaks}{\stackrel{*}{\rightharpoonup}}
\newcommand{\barw}{{\hspace*{.1em}\overline{\hspace*{-.1em}\W\hspace*{-.1em}}\hspace*{.1em}}}
\newcommand{\barf}{{\hspace*{.1em}\overline{\hspace*{-.1em}f}}}
\newcommand{\barz}{{\hspace*{.05em}\overline{\hspace*{-.05em}z}}}
\newcommand{\baru}{{\hspace*{.05em}\overline{\hspace*{-.05em}u}}}
\newcommand{\barF}{{\hspace*{.25em}\overline{\hspace*{-.25em}\mathcal{F}\hspace*{-.05em}}\hspace*{.05em}}}

\newcommand{\DeformationFile}[4][]{\includegraphics[#1]{Skew_MU#2TR#3_Def-#4}}
\newcommand{\PictureFile}[2][]{\includegraphics[#1]{#2}}

\begin{document}

\begin{center}
\Large\bfseries
Quasistatic normal-compliance contact problem of\\
visco-elastic bodies with Coulomb friction\\
implemented by QP and SGBEM.
\end{center}

\bigskip

\begin{center}
\bfseries
Roman Vodi\v cka\footnote{Technical University of Ko\v sice, Civil Engineering Faculty,
Vysoko\v skolsk\' a 4, 042 00 Ko\v sice, Slovakia, {\texttt{  roman.vodicka@tuke.sk}}},
Vladislav Manti\v{c}\footnote{University of Seville, School of Engineering,
Camino de los Descubrimientos s/n, 41092 Seville, Spain, {\texttt{  mantic@etsi.us.es}}},
Tom\'{a}\v{s} Roub\'\i\v{c}ek\footnote{Charles University, Mathematical Institute, Sokolovsk\'a 83,
CZ-186~75~Praha~8, Czech Republic, {\tt roubicek@karlin.mff.cuni.cz}}$^,$\footnote{Institute of Thermomechanics of the Czech Acad. Sci., Dolej\v skova~5,
CZ--182~08 Praha 8, Czech Rep.}
\end{center}

\begin{abstract}
The quasistatic normal-compliance contact problem of isotropic homogeneous linear visco-elastic bodies with Coulomb friction at small strains in Kelvin-Voigt rheology
is considered.
The discretization is made by a semi-implicit formula in time and the Symmetric Galerkin Boundary Element Method (SGBEM) in space, assuming that the ratio of the viscosity and   elasticity moduli is a given relaxation-time coefficient.
The obtained recursive minimization problem, formulated only on the contact boundary, has a nonsmooth cost function.
If the normal compliance responds linearly and the 2D problems are considered, then the cost function is piecewise-quadratic, which after a certain transformation gets the quadratic programming (QP) structure.
However,  it would lead to second-order cone programming in 3D problems.
Finally, several computational tests  are presented and analysed, with additional discussion on numerical stability and convergence of the involved approximated Poincar\'e-Steklov operators.

\medskip

\noindent{\textbf{Keywords}}:
contact mechanics, tribology, evolution variational inequalities, numerical approximation, boundary element method (BEM), mathematical programming, computational simulations.

\medskip

\noindent{\textbf{AMS Classification}}:
35Q90, 
49N10, 
65K15, 
65M38, 
74A55, 
74S15, 
90C20. 
\end{abstract}

\section{Introduction}

Mathematically justified modelling and efficient numerical solution of unilateral {\it contact problems} with friction is very challenging.
It is now well recognized that the combination of the {\it Coulomb friction} with unilateral Signorini contact is very difficult and essentially there are no satisfactory  mathematical results (i.e.\ not only for small coefficient of friction  or small data or short time) available in the literature.
On the other hand, in engineering a so-called {\it normal compliance} is  a well-accepted compromise (and sometimes considered even as a more realistic variant) to the Signorini problem.
This compliance concept is based on an idea of a rough contact surface which is intimately related with the concept of friction resulting microscopically from a certain asperity of the surface and which, when pressed by a normal force, exhibits certain elastic response, macroscopically demonstrated as a certain interpenetration, cf.~\cite[Sect.\,11.4.2]{Kikuchi}, \cite[Sect.\,4.10]{Maug00CARE},  \cite[Chap.\,5]{Wrig06CCM} or \cite{Tworzydlo1998}.
This normal-compliance concept is also well amenable for existence analysis, cf.~\cite{AndKla01rev,Klarbring1988,Rabier1986}, and for analysis of convergence of numerical approximation, although its algorithmic implementation might be more difficult because (in contrast to the Signorini problem) does not lead directly to a quadratic-programming or a 2nd-order-cone-programming problem.

When nonlinearities are restricted to the boundaries, small-strain concept is adopted, and the material in the bulk is linear and homogeneous (so that the fundamental solutions of the underlying operators are at disposal), the {\it Boundary Element Method} (BEM) can be an efficient and very accurate alternative to the {\it Finite Element Method }(FEM) usually employed in engineering computations, e.g.~\cite{Wrig06CCM,Laursen02,Kikuchi,kucera05,dostal10A1,Haslinger04A1,Haslinger06A1,Krause09}.
In this paper, we confine ourselves to the mentioned situation in the bulk, so that all nonlinearities arising from the normal-compliance model and the Coulomb friction will be indeed only on the boundary or interfaces between the bodies.
It should be noted that Coulomb friction can be reduced to the solution of so-called Tresca friction, see~\cite{DoHaKu02IFPM,Haslinger04A1} for which error estimates are available, see~\cite{Gwinner13}.
Also, it is worthy of  mention that we neglect inertia.

There exist several quite successful approaches for the solution of contact problems by BEM, e.g.~\cite{blazquez021,eck99A1,ChMaSt08MBEM,GMSS11AFEB,MaiSte05FEMB,PanManGG2013} and references therein.
In our present approach we develop a new model, which uses normal compliance contact with visco-elastic behaviour of the bulk but implemented so that the elasto-static boundary element (BE) approach is used in combination with the quadratic programming (QP).

The paper is organized as follows.
After formulation of the problem in Section~\ref{Sec_IDP}, we made the following steps towards devising an approximation and an efficient algorithm in Section~\ref{Sec_Implementation}:
\begin{itemize}
\item
Discretization made by a semi-implicit formula in time (Sect.~\ref{SecImplTime}).
\item
Implementation of BEM, assuming the ratio of the viscosity and elasticity moduli  is a given relaxation-time coefficient (Sect.~\ref{SecImplTime-BIE}).
\item
The changes in energy formulation related to implementation by the {\it Symmetric Galerkin Boundary Element Method (SGBEM)}  (Sect.~\ref{Sec_ImplSpace}).
\item
A transformation of the nonsmooth incremental problems to get the quadratic programming (QP) structure in 2D contact problems with linearly-responding normal compliance (Sect.~\ref{Sec_Minimization}).
\end{itemize}

Let us emphasize that the method is robust in the sense that it avoids any nonconvex minimization, and even the minimizers at each particular time level can be found in a finite number of steps, and has a guaranteed numerical stability (a-priori estimates) for arbitrary coefficient of friction  and large data and long times, as shown in Appendix~\ref{Sec_Theory}.
Most importantly, involving a viscous rheology and the normal compliance ensures that the continuous friction model has a unique solution, which was proved by Han et al\@.~\cite{HaShSo01VNAQ} (see also \cite{HanSo02QCPV,Shillor2004,CaFeKu09elviscon}) for the Kelvin-Voigt rheology, used also here.
It is worth mentioning that  uniqueness results for linear elastic solids are only available when assuming  some bounds for contact problem parameters \cite{Hild2009,Klarbring1988,Rabier1986}. 
In fact, any `enough dissipative' rheology leading to a parabolic-type equation (e.g.Jeffreys' or the so-called 4-parameter solid) will serve similarly but the less dissipative rheologies leading to a rather hyperbolic equation (like Maxwell's, Boltzmann's  or Burgers'), would not guarantee well-posedness of the model.
Sometimes the former type of rheology is  referred to as solid-type and the latter as fluid-type, cf.~\cite{PaMaRo14BAMI}.

Eventually, in Section~\ref{Sec_NumExample}, the algorithm exploiting the above constructed incremental QP (solved here by a conjugate gradient based method) is tested in numerical simulations for various types of contacts, including receding contact and conforming contact problems.

The paper also includes tree appendices.
In Appendix~\ref{App_SGBEM}, we give some details of SGBEM implementation for spatial discretization of a suitable system of {\it Boundary Integral Equations} (BIE), in  Appendix~\ref{Sec_PCC}, we present the block structure of the displacement-traction map for interface variables, in Appendix~\ref{Sec_Theory}, we sketch some theoretical aspects of the model concerning numerical stability and convergence of the numerical implementation.

\section{Contact problem with friction at small strains}\label{Sec_IDP}

\def\rmn{{_{\rm N}}}
\def\rmt{{_{\rm T}}}
\def\rmn{{\rm n}}
\def\rmt{{\rm t}}
\def\bbD{{\mathbb D}}
\def\bbC{{\mathbb C}}
\def\DT{\dot}
\def\R{{\mathbb R}}
\def\P{{p}}
\def\W{{w}}
\def\Z{{z}}
\def\DIR{{g}}

\begin{figure}[ht]
\centering
\PictureFile{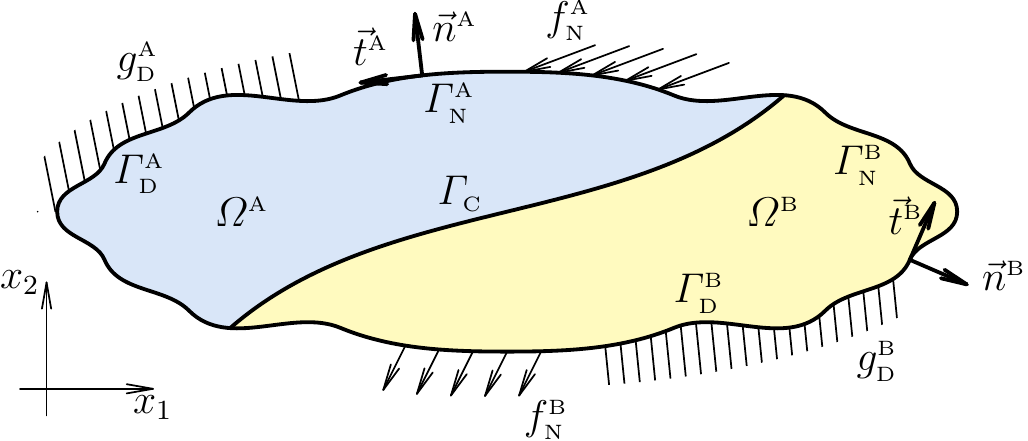}
\caption{Notation of two domains in contact and  various parts of their boundaries.}\label{Fig_Subdomains}
\end{figure}

For the notational simplicity, only two bodies in contact will be considered.
Let these bodies occupy domains $\Omega^\eta\!\subset\!\mathbb R^d$ $(\eta{=}{\rm A,B})$ with bounded Lipschitz boundaries $\partial\Omega^\eta{=}\Gamma^\eta$,
Figure~\ref{Fig_Subdomains}.
Later, we will consider only $d{=}2$ but some parts allow also $d{=}3$.
Let $\vec{n}^\eta$ denote the unit outward normal vector defined a.e.\ at $\Gamma^\eta$.
For $d{=}2$, $\Vector{t}^\eta$ denotes the unit tangential (surface) vector such that it defines anti-clockwise orientation of $\Gamma^\eta$.

The contact zone $\GC$ is defined as the common part of  $\GA$ and $\GB$, i.e. $\GC{=}\GA{\cap}\GB$.
Notice, that this definition, adopted for the sake of formulation simplicity, covers just receding and conforming contact types.
Nevertheless, the present computational implementation covers also the case of advancing contact.
The Dirichlet and Neumann boundary conditions are defined on the outer boundary parts, respectively, prescribing displacements as ${\DIR}\D^\eta={\DIR}\D^\eta(t)$ on $\GD^\eta$ and tractions as $f\N^\eta=f\N^\eta(t)$ on $\GN^\eta$ at a current time $t$, cf.~(\ref{BVP}b,c) below.
Of course, we assume that $\GD^\eta$ and $\GN^\eta$ are disjoint and relatively open parts of $\Gamma^\eta$, disjoint also with $\GC$, and
$\Gamma^\eta{=}\overline\Gamma_\Dir^\eta{\cup}\overline\Gamma_{_{\rm N}}^\eta{\cup}\overline\Gamma_{_{\rm C}}$.
We also assume that the Dirichlet part is far from the contact boundary, i.e. $\overline\Gamma_\Dir^\eta{\cap}\overline\Gamma_{_{\rm C}}{=}\emptyset$.
The difference, i.e.\ a gap, on $\GC$ of (traces of) functions defined on
$\OA$ and $\OB$ will be denoted by $\JUMP{\cdot}{}$.
In particular,  the gap of displacements on the contact boundary $\GC$ means $\JUMP{\Vector{u}}{}:=\Vector{u}\AA|_{\GC}^{}-\Vector{u}\BB|_{\GC}^{}$.
Further, we will use the gap of the normal and tangential displacements
\begin{equation}\label{Eq_NormalGap}
\JUMP{\Vector{u}}{\rmn}:=\JUMP{\Vector{u}}{}{\cdot}\vec{n}\BB
=-\JUMP{\Vector{u}}{}{\cdot}\vec{n}\AA
=-\Vector{u}\BB|_{\GC}^{}{\cdot}\vec{n}\BB-\Vector{u}\AA|_{\GC}^{}{\cdot}\vec{n}\AA, \qquad
\JUMP{\Vector{u}}{\rmt}:=\JUMP{\Vector{u}}{}-\JUMP{\Vector{u}}{\rmn}{\vec{n}}\BB
\end{equation}
We also use the convention that the `dot' will stand for the partial time derivative $\frac{\partial}{\partial t}$.

The classical formulation of the evolution boundary-value problem we will consider is the following:
\begin{subequations}\label{BVP}\begin{align}\label{BVP-equilibrium}
&\mathrm{div}\,\sigma^\eta=0,\ \ \sigma^\eta=
\bbD^\eta e(\DT{u}^\eta){+}\bbC^\eta e(u^\eta),\ \
e(u):=\mbox{$\frac12$}(\nabla u)^\top\!{+}\mbox{$\frac12$}\nabla u\!\!\!
&&\text{on }\Omega^\eta,\ \eta={\rm A,B},
\\\label{BVP-Dirichlet}
&\Vector{u}^\eta={\DIR}_\Dir^\eta(t)
&&\text{on }\GD^\eta,
\\\label{BVP-Neuman}
&
\P^\eta=\Vector{f}\N^{\eta}(t)
\ \ \ \text{ with the traction }\ \P^\eta:=\sigma^\eta\vec{n}^\eta
&&\text{on }\GN^\eta,
\\
&\!\!\!\!\left.\begin{array}{l}
\JUMP{\sigma}{}\vec{n}\BB=0\ \ \ \text{ and }\ \
\sigma_\rmn^{}-\gamma'\big(\JUMP{u}{\rmn}\big)=0\ \ \text{ with }\
\sigma_\rmn:=\P\BB{\cdot}\vec{n}\BB,
\\[.2em]
|\sigma_\rmt^{}|\le-\mu\sigma_\rmn^{}\ \ \ \text{ with }\
\sigma_\rmt^{}:=\P\BB-\sigma_\rmn^{}\vec{n}\BB,
\\[.2em]
|\sigma_\rmt^{}|<-\mu\sigma_\rmn^{}\ \ \ \Rightarrow\ \ \ \JUMP{\DT u}{\rmt}=0,
\\[.2em]
|\sigma_\rmt^{}|=-\mu\sigma_\rmn^{}\ \ \ \Rightarrow\ \ \ \exists\lambda\ge0:\ \
\sigma_\rmt^{}=\lambda\JUMP{\DT u}{\rmt},
\end{array}\right\}
&&\text{ on }\GC;
\label{eq6:friction-Robin-BC}
\end{align}
with $\bbC^\eta$ (resp.\ $\bbD^\eta$) the positive-definite symmetric 4th-order tensors of elastic (resp.\ viscous) moduli on $\eta$-domain, $\mu>0$ a coefficient of friction  and $\gamma:\R\to\R^+$ a decreasing smooth and convex function describing the compression energy of the normal-compliance contact.

The normal and tangential components of the traction stress defined in \eqref{BVP-Neuman} are denoted by $\sigma_\rmn^{}$ and $\sigma_\rmt^{}$ in \eqref{eq6:friction-Robin-BC}, respectively.
The evolution boundary-value problem \eqref{BVP} is to hold for a.a.\ time instances $t\in[0,T]$ with $T$ a fixed time horizon.
We will further consider an initial-boundary-value problem by prescribing the initial conditions at time $t{=}0$ for the  displacement:
\begin{align}\label{Eq_IC}
\Vector{u}^\eta(0,\cdot)=\Vector{u}_{0}^\eta\qquad\text{ on }\ \Omega^\eta,\ \ \eta={\rm A,B}.
\end{align}
\end{subequations}
It should be emphasized that the classical formulation \eqref{BVP} is rather
formal in some aspects and gain a sense only in a suitable weak
formulation.
In particular, the jump of stress vector $\JUMP{\sigma}{}\vec{n}\BB$ is well defined only in the space $H^{-1/2}(\GC;\R^d)$, and is imposed to equal zero,
in contrast to $\sigma_\rmn$ which has a good meaning in the space $L^{2-2/d}(\GC)$ if $\gamma'$ is qualified as \eqref{e4:friction-small-gamma} below,
being equal to $\gamma'(\JUMP{u}{\rmn})$.
It should be noted that the jump of the stress tensor $\JUMP{\sigma}{}$ is quite difficult to define and it is actually not needed in the formulation, additionally it does not have to vanish, while  the jump of the stress vector does.

As for the normal-compliance model, for some purposes of this paper, we can assume that, for some $2\le q<2(d-1)/(d-2)$ for $d{\ge}3$, otherwise $2{\le}q{<}\infty$, satisfying the following growth/continuity conditions:
\begin{subequations}\label{e4:friction-small-gamma}
\begin{align}\label{e4:friction-small-gamma1}
\exists\,
C\!\in\!\R\ \ &\forall\,u\!\in\!\R:\qquad\big|\gamma'(u)\big|\le C|u|^{q-1},\ \text{ and}
\\&\label{e4:friction-small-gamma2}
\forall\,u,\widetilde u\!\in\!\R:\quad\big|\gamma'(u){-}\gamma'(\widetilde u)\big|\le
C\big(|u|^{q-2}\!+|\widetilde u|^{q-2}\big)\big|u{-}\widetilde u\big|.
\end{align}\end{subequations}
Typically, $\gamma((0,+\infty])=0$ for unilateral contact and, e.g.,
$\gamma(u)=K|u^-|^q$ satisfies \eqref{e4:friction-small-gamma} provided
$q{\ge}2$. If $d{=}2$, $\gamma$ may have an arbitrary polynomial growth,
which may imitate a finite interpenetration. Let us remark that the finite
interpenetration up the depth, say, $u_{\rm f}{<}0$, means that
$\lim_{v{\searrow}u_{\rm f}}\gamma(v)=\infty$ and $\gamma(v)=\infty$ for
$v{\le}u_{\rm f}$. This is a quite realistic replacement of the Signorini
contact, although any rigorous analysis for the evolutionary model does not
seem to be  available yet. Anyhow, especially for an efficient implementation
by QP, we will later confine ourselves to a piecewise quadratic case,
i.e.\ in particular $q{=}2$.

For lucidity, it is noteworthy that the evolution boundary-value problem
\eqref{BVP-equilibrium} with \eqref{eq6:friction-Robin-BC}, when written in
a weak formulation, takes the form of a doubly-nonlinear evolution governed
by a convex stored energy $\Functional{E}$ and a nonsmooth convex potential
of dissipative forces $\Functional{R}(u;\cdot)$. In the so-called
Biot-equation form, it reads as the differential inclusion
\begin{equation}\label{Eq_InclusionVisco}
\partial_{\dot{\Vector{u}}}\Functional{R}(\Vector{u};\dot{\Vector{u}})
+\partial_{\Vector{u}}\Functional{E}({\DIR}_\Dir(t);\Vector{u})\ni
\Functional{F}(t),
\end{equation}
for $u:=(u\AA,u\BB)$ and ${\DIR}_\Dir=({\DIR}_\Dir\AA,{\DIR}_\Dir\BB)$, where the symbol $\partial$ refers to partial subdifferentials relying on the convexity of
$\Functional{R}(u;\cdot)$ and $\Functional{E}({\DIR}_\Dir;\cdot)$.
It involves the functionals $\Functional{E}$, $\Functional{R}$, and $\Functional{F}$ in the specific form:
\begin{subequations}\label{Eq_fun2}\begin{align}
 \label{Eq_Efun1}
\noindent
&\Functional{E}({\DIR}_\Dir;\Vector{u}):=\Functional{E}_{_\text{B}}({\DIR}_\Dir;\Vector{u})
+\Functional{E}\C(\Vector{u})\ \ \text{ with }\ \ \Functional{E}\C(\Vector{u}):=
\int_{\GC}\!\!\gamma\big(\JUMP{u}{\rmn}\big)\,\dd\Gamma,
\\ \label{Eq_Efun2}&
\Functional{E}_{_\text{B}}({\DIR}_\Dir;\Vector{u}):=\begin{cases}
\displaystyle{
\sum_{\eta={\rm A,B}}\int_{\Omega^\eta}\frac12 e(u^\eta){:}\bbC^\eta{:}e(u^\eta)\,\dd\Omega
}
&\text{if }\Vector{u}^\eta|_{\GD^\eta}^{}={\DIR}_\Dir^\eta,\\[-.3em]
\qquad\qquad\qquad\infty&\text{else},
\end{cases}
\\&\label{Eq_RfunFriction}
\Functional{R}(\Vector{u};\dot{\Vector{u}}):=
\Functional{R}_1^{}(\Vector{u};\dot{\Vector{u}})+\Functional{R}_2^{}(\dot{\Vector{u}})\ \ \text{ with }\ \
\Functional{R}_1^{}(\Vector{u};\dot{\Vector{u}}):=\int_{\GC}\!\!-\mu\gamma'\big(\JUMP{u}{\rmn}\big)\big|\JUMP{\DT u}{\rmt}\big|\,\dd\Gamma,
\\&\label{Eq_RfunVisco}
\Functional{R}_2^{}(\dot{\Vector{u}}):=\sum_{\eta={\rm A,B}}\int_{\Omega^{\eta}}
\frac12 e(\DT u^\eta){:}\bbD^\eta{:}e(\DT u^\eta)\,\dd\Omega,
\\&\label{Eq_Ffun2}
\big\langle\Functional{F}(t),v\big\rangle:=\sum_{\eta={\rm A,B}}
\int_{\GN^\eta}\!\!
\Vector{f}\N^{\eta}(t)
v^\eta\,\dd\Gamma\qquad
\forall v=(v\AA,v\BB),\ v^\eta\!\in\!H^1(\Omega^\eta;\R^d),\ \eta={\rm A,B}.
\end{align}\end{subequations}
Note that in \eqref{Eq_Efun1} and \eqref{Eq_Efun2} we split $\Functional{E}$ into the bulk part $\Functional{E}_{_\text{B}}$ and the contact-boundary part $\Functional{E}_{_\text{C}}$ to be used later especially in Section~\ref{SecImplTime-BIE}, while in \eqref{Eq_RfunFriction} and \eqref{Eq_RfunVisco} we denoted the split of the functional $\Functional{R}$ according to the degree of homogeneity of its terms with respect to the displacement rate (the subscript index).
The dissipation rate reflects these different homogeneity degrees, namely $\langle\Functional{R}(\Vector{u};\dot{\Vector{u}}),\dot{\Vector{u}}\rangle
=\Functional{R}_1^{}(\Vector{u};\dot{\Vector{u}})+2\Functional{R}_2^{}(\dot{\Vector{u}})$, as occurs in \eqref{Eq_Ebilance} below.

We intentionally do not want to write (\ref{BVP}a-d) in terms of some shifted solution satisfying the homogeneous Dirichlet conditions on $\GD^\eta$ because this (otherwise usual) procedure would give rise a non-homogeneous right-hand side in the transformed equilibrium equation \eqref{BVP-equilibrium} and, thus, it would not
be well suited with BEM for spatial discretization we want to use later.

The weak formulation of \eqref{Eq_InclusionVisco} is simple if $\DT{\DIR}_\Dir=0$, namely
\begin{align}\nonumber
&\forall_{\text{a.a.}}t\!\in\![0,T]\ \
\forall\,\wt u\!\in\!H^1(\Omega\AA{\cup}\Omega\BB;\R^d),\ \ \wt u|_{\GD}^{}=0:
\\&\qquad\qquad
\big\langle\partial_{\Vector{u}}\Functional{E}({\DIR}_\Dir;\Vector{u}(t))-\Functional{F}(t),
\wt u-\DT{\Vector{u}}(t)\big\rangle+\Functional{R}(\Vector{u}(t);\wt u)
\ge\Functional{R}(\Vector{u}(t);\dot{\Vector{u}}(t)).
\label{Eq_Inclusion-weak}
\end{align}
For a general case $\DT{\DIR}_\Dir\ne0$, like in \cite{ColVis90CDNE},
one should use the formulation as a system:
\begin{subequations}\label{Eq_Inclusion-weak+}\begin{align}
&\forall_{\text{a.a.}}t\!\in\![0,T]:\ \ \ \ \xi(t)+\zeta(t)=\Functional{F}(t),\ \ \text{ and }
\\&\forall\,\wt u\!\in\!H^1(\Omega\AA{\cup}\Omega\BB;\R^d):\ \ \ \big\langle
\xi(t),\wt u-\DT{\Vector{u}}(t)\big\rangle+\Functional{R}(\Vector{u}(t);\wt u)
\ge\Functional{R}(\Vector{u}(t);\dot{\Vector{u}}(t)),
\\
&\forall\,\wt u\!\in\!H^1(\Omega\AA{\cup}\Omega\BB;\R^d),
\ \wt u|_{\GD}^{}\!=0:\ \ \big\langle\zeta(t),\wt u{-}{\Vector{u}}(t)\big\rangle+
\Functional{E}({\DIR}_\Dir(t);\wt u)\ge\Functional{E}({\DIR}_\Dir(t);\Vector{u}(t)).
\end{align}\end{subequations}
However, one of the advantages of BIE and BEM is that the Dirichlet conditions are `hidden' into the stored-energy functional so that time-varying Dirichlet condition does not represent any problem if $\Gamma_\Dir$ is far from $\GC$, cf.~\eqref{Biot-disc-g-h-weak} and \eqref{Biot-disc-g-h-disc}.

Here, when we want  to write either the energy balance or \eqref{Eq_Inclusion-weak+} in a simpler form \eqref{Eq_Inclusion-weak}, we  use such shift only for testing \eqref{Eq_InclusionVisco}.
Let $u_\Dir$ define a (suitably smooth) prolongation of the boundary condition ${\DIR}_\Dir$.
Since $\overline\Gamma_\Dir^\eta{\cap}\overline\Gamma_{_{\rm C}}{=}\emptyset$ is assumed, we can always consider $\JUMP{u_\Dir}{}=0$ on $\GC$ so that this transformation
does not influence \eqref{eq6:friction-Robin-BC} and then also the $\Functional{R}_1^{}$-term in \eqref{Eq_InclusionVisco+} below.
Then, $u_{\rm z}:=u-u_\Dir$ has zero traces on $\GD$ and, in view of \eqref{Eq_InclusionVisco}, satisfies
\begin{equation}\label{Eq_InclusionVisco+}
\partial_{\dot{\Vector{u}}}\Functional{R}_1^{}(\Vector{u}_{\rm z};\DT u_{\rm z})
+\Functional{R}_2'(\DT u_{\rm z}{+}\DT u_\Dir)
+\partial_{\Vector{u}}\Functional{E}({\DIR}_\Dir(t);u_{\rm z}{+}u_\Dir)\ni
\Functional{F}(t),
\end{equation}
where $\Functional{R}_2'$ stands for G\^ateaux derivative.
Now, testing \eqref{Eq_InclusionVisco+} by $\dot{\Vector{u}}_{\rm z}$ having zero traces on $\GD$ becomes legal and renders
\begin{multline}\label{Eq_Eprebilance}
\sum_{\eta={\rm A,B}}\left(
\int_{\Omega^\eta}\!\!
e({\Vector{u}}_{\rm z}^\eta+u_\Dir^\eta){:}
\bbC^\eta{:}
e(\dot{\Vector{u}}_{\rm z}^\eta)\,\dd\Omega\right)
+\int_{\GC}\!\!
\gamma'\big(\JUMP{u}{\rmn}\big)\Big(\JUMP{\DT u}{\rmn}
-
\mu
\big|\JUMP{\DT u}{\rmt}\big|\Big)\,
\dd\Gamma
\\[-.3em]+\sum_{\eta={\rm A,B}}\bigg(
\int_{\Omega^\eta}\!\!
e(\dot{\Vector{u}}_{\rm z}^\eta+\DT u_\Dir^\eta)
{:}
\bbD^\eta{:}
e(\dot{\Vector{u}}_{\rm z}^\eta)\,\dd\Omega
 \bigg)
=
\sum_{\eta={\rm A,B}}\int_{\GN^\eta} \Vector{f}\N^{\eta}\cdot\dot{\Vector{u}}_{\rm z}^\eta\,\dd\Gamma,
\end{multline}
and by integrating it over the time interval $[0;T]$ with the initial conditions \eqref{Eq_IC} provides the following energy balance, by substituting ${\Vector{u}}_{\rm z}$ according to its definition,
\begin{multline}\label{Eq_Ebilance}
\Functional{E}({\DIR}_\Dir(T){;}\Vector{u}(T))
+\int_0^T\!\!\left(\Functional{R}_1^{}(\Vector{u};\dot{\Vector{u}})+2\Functional{R}_2^{}(\dot{\Vector{u}})\right)\dd t\\[-.1em]
=\Functional{E}({\DIR}_\Dir(0){;}\Vector{u}_0)+
\!\int_0^T\!\!\bigg(\big\langle\Functional{F}(t),\dot{\Vector{u}}{-}\dot u_\Dir\rangle
+\sum_{\eta={\rm A,B}}\int_{\Omega^\eta}\!\!\left(
\bbC^\eta e(u^\eta)+
\bbD^\eta e(\DT u^\eta)\right)
{:}
e(\dot u_\Dir^\eta)\dd\Omega\bigg)\dd t.
\end{multline}
This expresses {\it energy conservation} in terms of the original
displacement $u$: the sum of the energy stored in the system at time $T$ and
the energy dissipated during the time interval $[0;T]$ due to friction and
viscous processes equals the energy stored at the system at time $0$ plus
the work made by the external boundary loadings ${\DIR}_\Dir^\eta$ and $f\N^\eta$
during the same time  interval.

\section{Numerical approximation and computer implementation}\label{Sec_Implementation}

The numerical procedure devised for solving the above problem considers time and spatial discretizations separately, as usual.
The procedure is formulated in terms of the boundary data only, with the spatial discretization carried out by SGBEM.
We   restrict generality of the visco-elastic materials to homogeneous solids to allow for BEM-treatment, and further of the normal-compliance response  by assuming a quadratic one.
Eventually we also confine ourselves to 2D problems (i.e.\ $d=2$) in order to facilitate an efficient implementation by SGBEM combined with quadratic programming (QP); note that in 3D problems, the structure of the linear constraints \eqref{Eq_MoscoRestrictions} and \eqref{Eq_DiscreteRestrictions} below cannot be obtained, although efficient tools do exist for this case, too, cf.~Remark~\ref{rem-d=3} below.

\subsection{Semi-implicit discretization in time}\label{SecImplTime}

\def\TAU{\tau}

For the sake of simplicity, the time-stepping scheme is initially defined by a fixed time step size $\TAU$ such that $T/\TAU$ is integer.
Yet, it should be remarked that in specific calculations mixing stick/slip/jump regimes with typically very different time scale, a certain adaptivity is desirable and here very simple (as the inertial term is neglected).
The velocity is approximated by the finite difference $\dot{\Vector{u}}\approx(\Vector{u}^k{-}\Vector{u}^{k{-}1})/\TAU$, where $\Vector{u}^k$ denotes the approximate solution at the discrete time $k\TAU$.
We discretise the differential inclusion \eqref{Eq_InclusionVisco} by the {\it semi-implicit recursive scheme}
\begin{equation}\label{Biot-disc}
\partial_{\dot{\Vector{u}}}\Functional{R}\Big(\Vector{u}^{k{-}1};
\frac{\Vector{u}^k{-}\Vector{u}^{k{-}1}}{\TAU}\Big)
+\partial_{\Vector{u}}\Functional{E}({\DIR}_\Dir(k\TAU);\Vector{u}^k)\ni
\Functional{F}(k\TAU),\ \ \ \ \ k=1,...,T/\TAU,\ \ \ \ \ \Vector{u}^0=u_0.
\end{equation}
We consider $\TAU{>}0$ fixed thorough the whole section and thus, for notational simplicity, we will not explicitly emphasize it when writing $\Vector{u}^k$.
Occurrence of $\Functional{R}(\Vector{u}^{k{-}1};\cdot)$ rather than $\Functional{R}(\Vector{u}^k;\cdot)$, i.e.\ the semi-implicit formula rather than fully-implicit one, brings a {\it variational structure} to the incremental recursive problem \eqref{Biot-disc} in the sense that \eqref{Biot-disc} represents the first-order optimality condition for the strictly convex nonsmooth functional
\begin{equation}\label{Eq_DiscreteFunctional}
\Functional{H}^k(\Vector{u}):=
\TAU\Functional{R}\Big(\Vector{u}^{k{-}1};\frac{\Vector{u}{-}\Vector{u}^{k{-}1}}{\TAU}\Big)
+\Functional{E}({\DIR}_\Dir(k\TAU);\Vector{u})-\big\langle\Functional{F}(k\TAU),\Vector{u}\big\rangle.
\end{equation}
Thus, we can take the unique minimizer of $\Functional{H}^k$ as the solution to \eqref{Biot-disc}.

An analogy to~\eqref{Eq_Ebilance} can be obtained by making a test of~\eqref{Biot-disc} with the transformation by shifting the solution by $u_\Dir^k:=u_\Dir(k\TAU)$ into the form like \eqref{Eq_InclusionVisco+}.
Like the test of \eqref{Eq_InclusionVisco+} by $\DT u_{\rm z}$, we now make a test by ${\Vector{u}}_{\rm z}^k{-}{\Vector{u}}_{\rm z}^{k-1}$ at the time instant $k\TAU$, where ${\Vector{u}}_{\rm z}^k:=u^k-u_\Dir^k$.
In view of \eqref{Eq_fun2}, it yields:
\begin{align}\nonumber
&\int_{\GC}-\mu\gamma'\big(\JUMP{\Vector{u}_{}^{k-1}}{\rmn}\big){\cdot}\big|
\JUMP{{\Vector{u}}_{}^k{-}{\Vector{u}}_{}^{k-1}}{\rmt}\big|\,
\dd\Gamma
\\\nonumber
&\qquad+\sum_{\eta={\rm A,B}}\frac{1}{\TAU}\int_{\Omega^\eta}
e({\Vector{u}}_{\rm z}^{\eta,k}{-}{\Vector{u}}_{\rm z}^{\eta,k-1}+
u_\Dir^{\eta,k}{-}
u_\Dir^{\eta,k-1}){:}
\bbD^\eta{:}
e({\Vector{u}}_{\rm z}^{\eta,k}{-}{\Vector{u}}_{\rm z}^{\eta,k-1})\dd\Omega
\\\nonumber
&\qquad+\sum_{\eta={\rm A,B}}\int_{\Omega^\eta}
e({\Vector{u}}_{\rm z}^{\eta,k}+
u_\Dir^{\eta,k}){:}
\bbC^\eta{:}e({\Vector{u}}_{\rm z}^{\eta,k}{-}{\Vector{u}}_{\rm z}^{\eta,k-1})\dd\Omega\\&\qquad
+\int_{\GC}\gamma'\big(\JUMP{\Vector{u}_{}^k}{\rmn}\big){\cdot}\left(\JUMP{{\Vector{u}}_{}^k{-}{\Vector{u}}_{}^{k-1}}{\rmn}\right)\,\dd\Gamma
=\sum_{\eta={\rm A,B}}\int_{\GN^\eta}f\N^\eta{\cdot}({\Vector{u}}_{\rm z}^k{-}{\Vector{u}}_{\rm z}^{k-1})\,\dd\Gamma.\label{Eq_EbalanceDiscr}
\end{align}
It provides the following energy estimate, returning back to the original $\Vector{u}^k$:
\begin{align}\nonumber
&\Functional{R}_1^{}(\Vector{u}^{k-1};{\Vector{u}}^{k}{-}{\Vector{u}}^{k-1})
+\frac{2}{\TAU}\Functional{R}_2^{}({\Vector{u}}^{k}{-}{\Vector{u}}^{k-1})
+\Functional{E}({\DIR}_\Dir(k\TAU){,}\Vector{u}^k)
\leq\Functional{E}({\DIR}_\Dir((k{-}1)\TAU){,}\Vector{u}^{k-1})
\\\nonumber
&\qquad\qquad
+\sum_{\eta={\rm A,B}}\int_{\Omega^\eta}\left(\bbC^\eta e\left(\Vector{u}^{\eta,k-1}\right)
+\bbD^\eta e\left(\frac{{\Vector{u}}^{\eta,k}{-}{\Vector{u}}^{\eta,k-1}}{\TAU}\right)\right){:}e(u_\Dir^{\eta,k}{-}
u_\Dir^{\eta,k-1})\,\dd\Omega\\
&\qquad\qquad
+\sum_{\eta={\rm A,B}}\frac12\int_{\Omega^\eta}e(u_\Dir^{\eta,k}{-}u_\Dir^{\eta,k-1}){:}
\bbC^\eta{:}e(u_\Dir^{\eta,k}{-}u_\Dir^{\eta,k-1})\,\dd\Omega
+\big\langle\Functional{F}(k\TAU),\Vector{u}^{k}{-}{\Vector{u}}^{k-1}\big\rangle.
\label{Eq_EbilanceDiscr}
\end{align}
The above inequality is useful when we permit the variations of the time step parameter $\TAU$ between  subsequent time steps providing a form of adaptivity in the algorithm: the inequality is required to be satisfied as an approximation of equality with a given tolerance $\epsilon$ -- if the difference is greater than $\epsilon$, the time step is diminished, if the difference is much more less than $\epsilon$, the time step is enlarged, otherwise it is kept the same.

In the derivation of the inequality~\eqref{Eq_EbilanceDiscr}, the convexity of the functional $\Functional{E}(t,\cdot)$ from \eqref{Eq_Efun1} and \eqref{Eq_Efun2} has been used, which provides, e.g., the following relation for the pertinent boundary term:
\begin{equation}\label{Eq_ConvexityF}
\int_{\GC}\!\!\gamma'\big(\JUMP{\Vector{u}^k}{\rmn}\big){\cdot}
\JUMP{{\Vector{u}}^k{-}{\Vector{u}}^{k-1}}{\rmn}\,\dd\Gamma
\\
\geq
\int_{\GC}\!\!\gamma\big(\JUMP{\Vector{u}^k}{\rmn}\big)\,\dd\Gamma
-\int_{\GC}\!\!\gamma\big(\JUMP{\Vector{u}^{k-1}}{\rmn}\big)\,\dd\Gamma.
\end{equation}
The residuum in \eqref{Eq_EbilanceDiscr} can serve for varying adaptively the time step $\tau{>}0$ to keep this residuum under a pre-selected tollerance.
This, however, needs the energy conservation in the limit problem, cf.~Remark~\ref{rem-energy} in the Appendix~\ref{Sec_Theory}.

\subsection{Implementation of viscosity towards BIE}\label{SecImplTime-BIE}

Now we make a slight restriction of generality by assuming that
\begin{align}\label{visc-anzatz}
\bbD^\eta=\chi\bbC^\eta,\ \ \ \ \ \ \ \ \chi>0\ \text{given}.
\end{align}
In most materials (like metals) the relaxation time $\chi{>}0$ is small comparing to outer-loading time scale but we assume that it is rather large to affect the tangential contact with friction.
We assume $\chi$ to be the same for both subdomains; note that the potential \eqref{Eq_DiscreteFunctionalV} cannot be formulated if $\chi$ would be different in each subdomain; anyhow, some considerations (and computations) can be performed for non-potential problems, too, as suggested in \cite[Sect.3.3]{PaMaRo14BAMI}.
This allows for using the simple trick from \cite[Remark~6.2]{roubicek13A1}, which can be extended even for more complex rheologies as shown in \cite{PaMaRo14BAMI}, which facilitates usage of standard static BIEs even for quasi-static evolution problems.

The simple viscosity ansatz \eqref{visc-anzatz} is chosen in order to exploit the reformulation of the viscoelastic problem in the bulk in terms of an recursive static
elastic problem in the bulk, similarly as in \cite{PaMaRo14BAMI}, which is solved by the conventional elastostatic SGBEM here.
For this aim, let us introduce a new variable $\Vector{v}$ (a `fictitious' displacement) at the time level $k$, as
\begin{equation}\label{Eq_SubstV}
\Vector{v}^k:=\Vector{u}^k+\chi\frac{\Vector{u}^k-\Vector{u}^{k-1}}{\TAU}.
\end{equation}
Realizing that
\begin{equation}\label{Eq_SubstU}
\Vector{u}^k=\frac{\TAU}{\TAU{+}\chi}\Vector{v}^k+\frac{\chi}{\TAU{+}\chi}\Vector{u}^{k-1}
\ \ \ \ \text{ and also }\ \ \ \
\frac{\Vector{u}^k-\Vector{u}^{k-1}}\TAU=
\frac{\Vector{v}^k-\Vector{u}^{k-1}}{\TAU{+}\chi}
\end{equation}
and that, thanks to \eqref{visc-anzatz},
\begin{subequations}
\begin{align}
&\Functional{R}_2'\Big(\frac{\Vector{u}^k{-}\Vector{u}^{k-1}}\TAU\Big)
+\partial_{\Vector{u}}\Functional{E}_{_\text{B}}({\DIR}_\Dir(k\TAU);\Vector{u}^k)
=\partial_{\Vector{u}}\Functional{E}_{_\text{B}}(\wt {\DIR}_\Dir(k\TAU);\Vector{v}^k)
\\\label{Eq_AdmissibleV}
&\qquad\qquad\qquad\qquad\qquad\text{with }\ \wt {\DIR}_\Dir(k\TAU):=
{\DIR}_\Dir(k\TAU)+\frac{\chi}{\TAU}\bigl({\DIR}_\Dir(k\TAU)-{\DIR}_\Dir((k{-}1)\TAU)\bigr),
\end{align}\end{subequations}
we can re-write \eqref{Biot-disc} in terms of this new variable as:
\begin{multline}\label{Biot-disc-v}
\partial_{\dot{\Vector{u}}}\Functional{R}_1^{}\big(\Vector{u}^{k{-}1};
\Vector{v}^k{-}\Vector{u}^{k-1}\big)
+\partial_{\Vector{u}}\Functional{E}_{_\text{B}}(\wt {\DIR}_\Dir(k\TAU);\Vector{v}^k)\\
+\Functional{E}_{_\text{C}}'\Big(\frac{\TAU\Vector{v}^k{+}\chi\Vector{u}^{k-1}}{\TAU+\chi}
\Big)
\ni
\Functional{F}(k\TAU),\ \ \ \ k=1,...,T/\TAU;
\end{multline}
here, we have also used that $\Functional{R}_1^{}(\Vector{u}^{k{-}1};\cdot)$ is (positive) homogeneous of degree 1 so that $\partial_{\dot{\Vector{u}}}\Functional{R}_1^{}(\Vector{u}^{k{-}1};\cdot)$ is homogeneous of degree 0 and the factor $\frac1{\TAU+\chi}$ from \eqref{Eq_SubstU} disappears from the $\Functional{R}_1^{}$-term in \eqref{Biot-disc-v}.
The recursive scheme then involves also the first equality in \eqref{Eq_SubstU} and starts with $u^0=u_0$ as in \eqref{Biot-disc}.
It is important that all bulk integrals are now contained in the `static' part $\Functional{E}$.
The inclusion \eqref{Biot-disc-v} has again the variational structure employing
the convex, nonsmooth functional:
\begin{align}\nonumber
\Functional{K}^k(\Vector{v})&=
\sum_{\eta={\rm A,B}}\int_{\Omega^\eta}\frac12e(v^\eta){:}\bbC^\eta{:}e(v^\eta)\,\dd\Omega
-\sum_{\eta={\rm A,B}}\int_{\GN^\eta}f\N^\eta{\cdot}\Vector{v}^\eta\dd\Gamma
\\&\ \ \ \
+\int_{\GC}\!\!\bigg(\frac{\TAU{+}\chi}{\TAU}
\gamma\Big(\JUMP{\frac{\TAU\Vector{v}{+}\chi\Vector{u}^{k-1}}{\TAU+\chi}}{\rmn}\Big)
-
\mu\gamma'\big(\JUMP{\Vector{u}^{k-1}}{\rmn}\big)
\big|\JUMP{\Vector{v}-\Vector{u}^{k-1}}{\rmt}\big|\bigg)\,\dd\Gamma
\label{Eq_DiscreteFunctionalV}
\end{align}
where $\Vector{v}=(v\AA,v\BB)$.
Like before, \eqref{Biot-disc-v} represents the 1st-order optimality condition for $\Vector{v}^k$ to be a minimizer of $\Functional{K}^k(\cdot)$ subject to the boundary condition $\Vector{v}^{\eta,k}={\wt{\DIR}}_\Dir^\eta(k\TAU)$ for $\eta={\rm A,B}$.

The BIE-method used here is based on a formulation of the problem only on the
contact boundary $\GC$. To this goal, we abbreviate the displacement and the
fictitious displacement gaps by
\begin{subequations}\begin{align}
&&&&&\Vector{\Z}:=\JUMP{\Vector{u}}{},&&
\Vector{\Z}_\rmn^{}:=\JUMP{\Vector{u}}{\rmn},&&
\Vector{\Z}_\rmt^{}:=\JUMP{\Vector{u}}{\rmt},\quad\text{ and }&&&&
\\&&&&&\Vector{\W}:=\JUMP{\Vector{v}}{},&&
\Vector{\W}_\rmn^{}:=\JUMP{\Vector{v}}{\rmn},&&
\Vector{\W}_\rmt^{}:=\JUMP{\Vector{v}}{\rmt}.&&&&
\end{align}\end{subequations}
Moreover, we define the linear mapping $S:({\Vector{\DIR}}_\Dir,f\N,\W)\mapsto v$ with $v$ being the weak solution to the {\it Transmission Boundary Value Problem} (TBVP):
\begin{subequations}\label{BVP-k}\begin{align}
&&&{\rm div}(\bbC^\eta e(v^\eta))=0&&\text{in $\Omega^\eta$,\ \ $\eta={\rm A,B}$},&&
\\&&&v={{\Vector{\DIR}}}^\eta_\Dir
&&\text{on $\Gamma^\eta\D$,\ \ $\eta={\rm A,B}$},&&
\\&&&\P^\eta=f\N^\eta
\ \ \ \text{ with }\ \ \P^\eta=(\bbC^\eta e(v)){\cdot}\vec{n}^\eta
\!\!\!\!\!\!\!\!
&&\text{on $\Gamma^\eta\N$,\ \ $\eta={\rm A,B}$},&&
\\&&& \!\!\!\!\left.\begin{array}{l} \JUMP{v}{}=\W\\[.2em]
\P\AA+\P\BB=0\end{array}\right\}&&\text{on $\GC$}.&&
\end{align}
\end{subequations}
Later in Section~\ref{Sec_ImplSpace}-\ref{Sec_Minimization}, we will solve \eqref{BVP-k} approximately by SGBEM, see Appendix~\ref{App_SGBEM}, using the formulation for multi-domain problems in~\cite{vodicka07A1}.
To this goal, it is desirable to replace all bulk integrals in \eqref{Biot-disc-v}.
The resulted recursive inclusion will be then formulated only on $\GC$, being written exclusively in terms of $\Z$ and $\W$ only:
\begin{subequations}\label{Biot-disc-g-h}\begin{align}\label{Biot-disc-g-h+}
&\partial_{\DT{\W}}\mathcal{R}_1^{}
\big(\Z_\rmn^{k{-}1};\W_\rmt^k{-}\Z_\rmt^{k-1}\big)
+
\big[\mathcal{E}_{_\text{B}}\big]_{\W}'(\wt {\DIR}_\Dir(k\TAU);\W^k)
+\mathcal{E}_{_\text{C}}'\Big(\frac{\TAU\W_\rmn^k{+}\chi\Z_\rmn^{k-1}}{\TAU+\chi}
\Big)\ni
\mathcal{F}_\W'(k\TAU),
\\&\label{Biot-disc-g-h++}
\Z^k=\frac\TAU{\TAU{+}\chi}\W^k+\frac\chi{\TAU{+}\chi}\Z^{k-1},
\qquad\qquad k=1,...,T/\TAU;
\end{align}\end{subequations}
for the recursion \eqref{Biot-disc-g-h++}, cf.~\eqref{Eq_SubstU}.
Of course,  for $k{=}1$, we start with the initial condition $\Z^0{=}\JUMP{u_0}{}$, calculate $\W^1$ from \eqref{Biot-disc-g-h+}, then $\Z^1$ from \eqref{Biot-disc-g-h++}, and continue for $k{=}2$ etc..
In \eqref{Biot-disc-g-h+}, we have used
\begin{subequations}\label{def-of-R-E-E-F}\begin{align}
&\mathcal{R}_1^{}(\Z_{\rmn}^{};\DT{\W}_{\rmt}^{})
=\mathcal{R}_1^{}\big(\JUMP{u}{\rmn};\JUMP{\DT v}{\rmt}\big)
:=\Functional{R}_1^{}(u,\DT v),
&&\mathcal{E}_{_\text{B}}(\wt {\DIR}_\Dir;\W)
:=\Functional{E}_{_\text{B}}(\wt{\DIR}_\Dir;S(\wt {\DIR}_\Dir,f\N,\W)),
\\&\mathcal{E}_{_\text{C}}(\W_{\rmn})=\mathcal{E}_{_\text{C}}(\JUMP{v}{\rmn})
:=\Functional{E}_{_\text{C}}(v),
&&
\mathcal{F}(t,\W):=
\big\langle\Functional{F}(t),S(\wt{\DIR}_\Dir,f\N,\W)\big\rangle.
\end{align}\end{subequations}
Note that $\Functional{R}_1^{}$ and $\Functional{E}_{_\text{C}}$ depend only on the values on $\GC$ and thus the mapping $S$ is not used in the definition of $\mathcal{R}_1^{}$ and $\mathcal{E}_{_\text{C}}$.
Note also that, in contrast to $\Functional{E}_{_\text{B}}(\wt {\DIR}_\Dir;\cdot)$, the functional $\mathcal{E}_{_\text{B}}(\wt {\DIR}_\Dir;\cdot)$ is smooth.
The functional~\eqref{Eq_DiscreteFunctionalV} can be written as a function of $\Vector{\W}$ as
\begin{align}\nonumber
\mathcal{K}^k(\Vector{\W})&:=
\int_{\GC}\!\!\bigg(\frac{\TAU{+}\chi}{\TAU}
\gamma\Big(\frac{\TAU}{\TAU{+}\chi}\W_{\rm n}^{}
+\frac{\chi}{\TAU{+}\chi}\Z_{\rm n}^{k-1}\Big)
-
\mu\gamma'\big(
\Z_{\rm n}^{k-1}\big)
\big|\W_{\rmt}-\Z_{\rmt}^{k-1}\big|
\\&\nonumber\ \ \ \
-\frac12\Vector{\P}\BB
\left({\wt{\Vector{\DIR}}_\Dir}(k\TAU),{\Vector{f}\N}(k\TAU),\Vector{\W}\right)
\cdot\Vector{\W}\bigg)\,\dd\Gamma
-\sum_{\eta={\rm A,B}}\int_{\GN^\eta}\frac12\Vector{f}\N^\eta\cdot\Vector{v}^\eta
\left({\wt{\Vector{\DIR}}}_\Dir^{}(k\TAU),{\Vector{f}\N}(k\TAU),\Vector{\W}\right)
\dd \Gamma
\\&\ \ \ \
+\sum_{\eta={\rm A,B}}\frac12\int_{\Gamma^\eta\D}
\Vector{\P}^\eta\left({\wt{\Vector{\DIR}}}_\Dir(k\TAU),{\Vector{f}\N}(k\TAU),
\Vector{\W}\right) \cdot{\wt{\Vector{\DIR}}}^\eta_\Dir(k\TAU)\,\dd\Gamma
\label{Eq_ChangeOmegaToGamma}
\end{align}
where $\Vector{v}^\eta \left({\wt{\Vector{\DIR}}}_\Dir^{}(k\TAU),{\Vector{f}\N}(k\TAU),\Vector{\W}\right)$ and $\Vector{\P}^\eta\left({\wt{\Vector{\DIR}}}_\Dir(k\TAU),{\Vector{f}\N}(k\TAU),\Vector{\W}\right)$
denote the displacement and traction, respectively, resulting in TBVP \eqref{BVP-k} through the (generalized) Poincar\'e-Steklov operator $\calPt$ which is described in Appendix~\ref{Sec_PCC} and whose representation is given by the relation \eqref{Eq_SblockTwo}.

In fact, the expressions $\int_{\GC}\P{\cdot}\W\,\dd\Gamma$ and $\int_{\GD}\P{\cdot}{\wt{\Vector{\DIR}}}^\eta_\Dir\,\dd\Gamma$ are to be understood rather as $(H^{-1/2},H^{1/2})$-dualities than Lebesgue integrals.
We dare to keep this non-precise but more conventional notation in engineering throughout this article.

Let us further realize that $\mathcal{K}^k(\Vector{\W})=\Functional{K}^k(v)$ holds if $v$ solves the boundary-value problem \eqref{BVP-k}.
The optimality conditions for minimization of $\mathcal{K}^k$ represent a weak formulation of the discrete problem \eqref{Biot-disc-v}.
The advantage of \eqref{Biot-disc-g-h+} in contrast to \eqref{Biot-disc-v} is that the constraints coming from the Dirichlet conditions (which prevent to write the weak formulation of the original problem generally as in~\eqref{Eq_Inclusion-weak}) are eliminated.
To be more specific:
\begin{align}\nonumber
&\forall\,\wt{\W}\!\in\!H^{1/2}(\GC;\R^d):\ \
\mathcal{R}_1^{}\big(\Z_\rmn^{k{-}1};\wt{\W}_\rmt{-}\Z_\rmt^{k-1}\big)
+\Big\langle\mathcal{E}_{_\text{C}}'\Big(
\frac{\TAU\W_\rmn^k{+}\chi\Z_\rmn^{k-1}}{\TAU+\chi}\Big),
\wt{\W}_\rmn^{}{-}\W_\rmn^k\Big\rangle
\\[-.2em]&\hspace*{10em}
+\big\langle\P^k,\wt{\W}{-}\W^k\big\rangle
\ge\mathcal{R}_1^{}\big(\Z_\rmn^{k{-}1};\W_\rmt^k{-}\Z_\rmt^{k-1}\big)
+\big\langle\mathcal{F}_\W'(k\TAU),\wt{\W}-\W^k\big\rangle\label{Biot-disc-g-h-weak}
\end{align}
with the traction $\P^k$ being $\Vector{\P}^\eta\left({\wt{\Vector{\DIR}}}_\Dir(k\TAU),{\Vector{f}\N}(k\TAU),\Vector{\W}\right)$ from \eqref{Eq_ChangeOmegaToGamma}.

\subsection{Spatial discretization and SGBEM}\label{Sec_ImplSpace}

The role of the SGBEM in the present computational procedure is to provide a complete boundary-value solution from the given boundary data to calculate the elastic strain
energy in these domains.
Thus, at each time step and at each iteration of the minimization algorithm, the SGBEM code calculates, in a similar way as in~\cite{vodicka07A1},  unknown tractions $\Vector{\P}$ along $\GC\cup\GD$ and unknown displacements \Vector{v} along $\GC\cup\GN$,  assuming the displacement gap $\Vector{\W}$ on $\GC$ to be known from the used minimization procedure.

Here, we rely on the assumption that the material is homogeneous in each solid to allow for BEM-treatment.
Then, the present implementation of SGBEM, deduced from~\cite{vodicka07A1,vodicka11A1}, is briefly described in Appendix~\ref{App_SGBEM}.
In fact, the necessary SGBEM calculations are besad on the solution of~\eqref{RepresentationMatrix}.
The variables of the (modified)  displacement \Vector{v}, the displacement gap \Vector{\W} and the traction \Vector{\P}  appearing in~\eqref{Eq_ChangeOmegaToGamma} are approximated as given by the formula~\eqref{Eq_UTapproximation}, i.e.
\begin{equation*} 
\Vector{v}^\eta(x)=\sum_n\Vector{N}^\eta_{\psi n}(x)\VectorMatrix{v}^\eta_n,\quad
\Vector{\W}(x)=\sum_{n\C^{}}\Vector{N}\BB_{\psi {n\C^{}}}(x)\VectorMatrix{\W}_{n\C^{}},\quad
\Vector{\P}^\eta(x)=\sum_m\Vector{N}^\eta_{\varphi m}(x)\VectorMatrix{\P}^\eta_m, \eqno{~\eqref{Eq_UTapproximation}}
\end{equation*}
with nodal variables gathered in respective vectors $\VectorMatrix{v}^\eta_\ell$, $\VectorMatrix{\W}_{\ell}$ and $\VectorMatrix{\P}^\eta_\ell$.

The SGBEM-approximation of $\mathcal{K}^k$ is the discretized energy from the equation~\eqref{Eq_ChangeOmegaToGamma} with approximations~\eqref{Eq_UTapproximation},
let us denote it by $\mathcal{K}^k_h$ with $h>0$ referring to the mesh-size of the BE discretization of $\GC$, is then
\begin{subequations}\label{def-of-K-h},
\begin{align}\nonumber
\mathcal{K}^k_h(\VectorMatrix{\W})&:=
\int_{\GC}\!\!\bigg(\frac{\TAU{+}\chi}{\TAU}
\gamma\Big(\frac{\TAU}{\TAU{+}\chi}
[\Vector{\W}_h^{}]_\rmn^{}
+
\frac{\chi}{\TAU{+}\chi}
[\Vector{\W}_h^{k-1}]_\rmn^{}
\Big)
-
\mu\gamma'\big(
\VectorMatrix{\Z}_{\rmn,n\C^{}}^{k-1}\big)
\Big|
[\Vector{\W}_h^{}]_\rmt^{}-
[\Vector{\Z}_h^{k-1}]_\rmn^{}\Big|
\\&\nonumber\qquad\qquad\ -\frac12
\P_h\BB(\wt {\DIR}_\Dir(k\TAU),f\N(k\TAU),\W_h^{})
\cdot
\Vector{\W}_h^{}
\bigg)\,\dd\Gamma
\\&\nonumber\ \ \ \
-\sum_{\eta={\rm A,B}}\int_{\GN^\eta}\frac12\Vector{f}\N^\eta\cdot\Vector{v}_h^\eta
\left({\wt{\Vector{\DIR}}}_\Dir(k\TAU),{\Vector{f}\N}(k\TAU),\Vector{\W}_h^{}\right)\,\dd \Gamma
\\&\ \ \ \
+\sum_{\eta={\rm A,B}}\frac12\int_{\Gamma^\eta\D}
\P_h^\eta(\wt {\DIR}_\Dir(k\TAU),f\N(k\TAU),\W_h^{})
\cdot{\wt{\Vector{\DIR}}}^\eta_\Dir(k\TAU)\,\dd\Gamma,
\\&\text{with $\Vector{\W}_h=\Vector{\W}_h(\VectorMatrix{\W})$ from
\eqref{Eq_UTapproximation}, i.e.\
$\Vector{\W}_h(x)=
\sum_{n\C^{}}\Vector{N}\BB_{\psi {n\C^{}}}(x)\VectorMatrix{\W}_{n\C^{}}$,}
\end{align}\end{subequations}
and  with $\Vector{v}_h^\eta \left({\wt{\Vector{\DIR}}}_\Dir^{}(k\TAU),{\Vector{f}\N}(k\TAU),\Vector{\W}_h\right)$ and $\Vector{\P}_h^\eta\left({\wt{\Vector{\DIR}}}_\Dir(k\TAU),{\Vector{f}\N}(k\TAU),\Vector{\W}_h\right)$ obtained by the approximate (generalized) Poincar\'e-Steklov operator $\calPt_h$ of Appendix \ref{Sec_PCC}, also used in solution of contact problems by BEM~\cite{MaiSte05FEMB,MaiSte05hpBEM}.
Here also, $\Gamma\BB$ is  the master surface for displacements, similarly  as in~\cite{vodicka07A1}, and $\Vector{\W}_h^{k-1}\in V_h$ verifies the recursive formula  analogous to \eqref{Biot-disc-g-h++}, i.e.\ $\Vector{\Z}_h^{k}=\frac{\tau}{\tau+\chi}\Vector{\W}_h^{k}
+\frac{\chi}{\tau+\chi}\Vector{\Z}_h^{k-1}$.
The BEM-version of \eqref{Biot-disc-g-h-weak} seeks $w_h^k\in V_h$ (denoted here for notational simplicity without $h$) such that
\begin{subequations}\label{Biot-disc-g-h-disc}
\begin{multline}
\forall\,\wt{\W}\!\in\!V_h:\ \ \ \
\mathcal{R}_1^{}\big(\Z_\rmn^{k{-}1};\wt{\W}_\rmt{-}\Z_\rmt^{k-1}\big)
+\Big\langle
\mathcal{E}_{_\text{C}}'\Big(
\frac{\TAU\W_\rmn^k{+}\chi\Z_\rmn^{k-1}}{\TAU+\chi}\Big),
\wt{\W}_\rmn{-}\W_\rmn^k\Big\rangle\\
+\big\langle\P^k,\wt{\W}{-}\W^k\big\rangle
\ge\mathcal{R}_1^{}\big(\Z_\rmn^{k{-}1};\W_\rmt^k{-}\Z_\rmt^{k-1}\big)
+\big\langle\mathcal{F}_\W'(k\TAU),\wt{\W}-\W^k\big\rangle
\end{multline}
where $V_h=\{\W\in W^{1,\infty}(\GC;\R^d);\ \W\text{ in the form \eqref{Eq_UTapproximation}}\}\subset H^{1/2}(\GC;\R^d)$ and $\P_h^k$ is defined by $\Vector{\P}_h^\eta\left({\wt{\Vector{\DIR}}}_\Dir(k\TAU),{\Vector{f}\N}(k\TAU),\Vector{\W}_h\right)$ from \eqref{def-of-K-h}.
Of course, again we use the recursion \eqref{Biot-disc-g-h++}, which is written in the terms of the spatial discretizations as
\begin{align}
\VectorMatrix{\Z}^k=\frac\TAU{\TAU+\chi}\VectorMatrix{\W}^k
+\frac\chi{\TAU+\chi}\VectorMatrix{\Z}^{k-1},
\end{align}
\end{subequations}
provided the initial condition $u_0$ is assumed so that $\JUMP{u_0}{}\in V_h$ for $h>0$ sufficiently small (or it is approximated  by an element of the space $V_h$).

\subsection{QP-minimization algorithm in 2-dimensional problems}\label{Sec_Minimization}

Once all the boundary data are known, the energy of the state given by $\mathcal{K}^k$ in~\eqref{Eq_ChangeOmegaToGamma} can be calculated.
It is worth seeing how it can be carried out in the present implementation provided we confine ourselves to the most common case of linearly-responding compliance on $\GC$, i.e.\ piecewise quadratic $\gamma$:
\begin{align}\label{quadratic-compliance}
\gamma(g)=\frac{k_g}2\min(0,g)^2,\ \ \ \ k_g>0\text{ given}.
\end{align}
Further, except Remark~\ref{rem-d=3}, we confine ourselves to 2D problems where $\Functional{R}_1^{}(u,\cdot)$ has, after discretization, a polyhedral graph, and similarly also $\Functional{E}\C$,  with a help of \eqref{quadratic-compliance}, leads after discretization to polyhedral constraints, see \eqref{Eq_DiscreteRestrictions} below.
All this will allow for usage of efficient quadratic-programming (QP) algorithms.

First, let us re-consider the non-quadratic terms in $\mathcal{K}^k$.
A rather classical trick of handling such unpleasant terms by displacing them into a position of constraints by additional variables; \emph{Mosco}-type \emph{transformation}, although the original work \cite{Mosc67RTFB} was designed rather for the optimality conditions as variational inequalities, which works even for nonpotential problems.
For the rate-independent problems where the unpleasant terms are piecewise
affine, it was used in \cite[Lemma~4]{Roub02EMMP} and then also in
\cite[Lemma~4]{RouKru04MEMM}.
Here we have also piecewise quadratic terms, but they need only a simple modification to consider their square root to obtain again the linear-constraint structure but then add the square power of these auxiliary variables into the cost function.
Let these additional auxiliary variables handling the piecewise linear and the piecewise quadratic terms be denoted as $\alpha$ and $\beta$, respectively, and the following constraints hold for $\Vector{\W}=(\Vector{\W}_\rmt,\Vector{\W}_\rmn)$ and $(\alpha,\beta)$:
\begin{subequations}\label{Eq_MoscoRestrictions}
\begin{align}
&&&&\alpha-\Vector{\W}_\rmt&\geq-\JUMP{\Vector{u}^{k-1}}{\rmt},&
\alpha+\Vector{\W}_\rmt&\geq\JUMP{\Vector{u}^{k-1}}{\rmt},&&&&\\
&&&&\beta&\geq0,&
\beta+\Vector{\W}_\rmn&\geq-\frac{\chi}{\TAU}\JUMP{\Vector{u}^{k-1}}{\rmn}.&&&&
\end{align}
\end{subequations}
For the discretization, the approximation formulas for both auxiliary parameters $\alpha$ and $\beta$ given by  the pertinent BE mesh should be considered.
In what follows, the same mesh and approximation as used in~\eqref{Eq_UTapproximation} for displacements on the boundary part $\GC\BB$ is considered.
The approximation formulas can be written in the form
\begin{equation}\label{Eq_ABapproximation}
\alpha(x)=\sum_{n\C^{}}^{} N\BB_{\psi n\C^{}}(x)\alpha_{n\C^{}}^{},\qquad
\beta(x)=\sum_{n\C^{}}^{} N\BB_{\psi n\C^{}}(x)\beta_{n\C^{}}^{},
\end{equation}
where  $\alpha_{n\C^{}}^{}$, $\beta_{n\C^{}}^{}$ are the nodal unknowns associated to the node $x_{n\C^{}}\BB$.

The restrictions~\eqref{Eq_MoscoRestrictions} written for the discretization can be obtained as relations between the nodal values, realizing that the right-hand side terms $\wt{\Vector{u}}^{k-1}{=}\JUMP{\Vector{u}^{k-1}}{}$ can be expressed in terms of $\Vector{h}^{k-1}$ applying to them the transformation~\eqref{Eq_SubstU}.
If $n$ runs through all nodes in the possible contact zone $\GC$ with respect to the boundary $\GB$, the nodal relations read
\begin{subequations}\label{Eq_DiscreteRestrictions}
\begin{align}
&&&&\alpha_{n\C}^{}-\W_{\rmt,n\C^{}}&\geq -\wt{u}^{k-1}_{\rmt,n\C^{}},
&
\alpha_{n\C}^{}+\W_{\rmt,n\C^{}}&\ge\JUMP{\wt{u}^{k-1}_{n\C^{}}}{\rmt},
\\&&&&
\beta_{n\C}^{}&\geq0,&
\beta_{n\C}^{}+\W_{\rmn,n\C^{}}&\geq -\frac{\chi}{\TAU}\JUMP{\wt{u}^{k-1}_{n\C^{}}}{\rmn},&&&&
\end{align}
\end{subequations}
where in the nodal values at the node $n$ subscript indices $\rmt$ and $\rmn$ are used to distinguish between tangential and normal components, respectively.

Then, in terms of these new variables, ${\mathcal{K}}^k_h$ from \eqref{def-of-K-h} takes the following form  denoted   by $\widehat{\mathcal{K}}^k_h$:
\begin{multline}\label{Eq_DiscreteEnergyV}
\widehat{\mathcal{K}}^k_h(\VectorMatrix{\W},\boldsymbol\alpha,\boldsymbol\beta)
=\int_{\GC}\!\!\bigg(\frac12\frac{\TAU k_g}{\TAU{+}\chi}\Big(\sum_{n\C^{}} N\BB_{\psi,n\C^{}}(x)\beta_{n\C^{}}\Big)^2\\
\qquad+\mu k_g\Big(\sum_{n\C^{}} N\BB_{\psi,n\C^{}}(x)\beta^{k-1}_{n\C^{}}\Big)\Big(\sum_{n\C^{}} N\BB_{\psi,n\C^{}}(x)\alpha_{n\C^{}}\Big)\\
\qquad\qquad\ -\frac12\sum_{n\C^{}}\Vector{N}\BB_{\psi n\C^{}}(x)\VectorMatrix{\W}_{n\C^{}}\cdot\sum_r\Vector{N}\BB_{\varphi r}(x)\VectorMatrix{\P}\BB_r(\wt{\VectorMatrix{\DIR}}_\Dir,\VectorMatrix{f}\N,\VectorMatrix{\W})\bigg)\,\dd\Gamma\\
\qquad\ +\sum_{\eta={\rm A,B}}
\int_{\GA\D}\frac12\sum_{n\D^{}}\Vector{N}^\eta_{\psi,n\D^{}}(x)\wt{\VectorMatrix{\W}}^{\eta,k}_{n\D^{}}\cdot\sum_{n\D^{}}\Vector{N}^\eta_{\varphi,n\D^{}}(x)\VectorMatrix{\P}^{\eta, k}_{n\D^{}}(\wt{\VectorMatrix{\DIR}}_\Dir,\VectorMatrix{f}\N,\VectorMatrix{\W})\,\dd\Gamma\\
-\sum_{\eta={\rm A,B}}\int_{\GN^\eta}\frac12 \sum_{n\N^{}}\Vector{N}^\eta_{\psi,n\N^{}}(x)\VectorMatrix{v}^{\eta, k}_{n\N^{}}(\wt{\VectorMatrix{\DIR}}_\Dir,\VectorMatrix{f}\N,\VectorMatrix{\W})\cdot
\sum_{n\N^{}}\Vector{N}^\eta_{\varphi,n\N^{}}(x){\VectorMatrix{f}\N^{\eta, k}}_{n\N^{}}\,\dd\Gamma.
\end{multline}
where all the unknown data written as a function of $\Vector{\W}$ are obtained by the SGBEM code.

In minimization of the functional~\eqref{Eq_DiscreteEnergyV} with the restrictions~\eqref{Eq_DiscreteRestrictions}, it may be useful to reformulate the problem in such a way that the restrictions change to bound constraints.
The left-column and right-column restrictions in~\eqref{Eq_DiscreteRestrictions} provide linearly independent constraints which can respectively be written in a matrix form as
\begin{equation}\label{Eq_SubstDef}
\begin{pmatrix} \Matrix{I}  & \Matrix{I} \\ \Matrix{I}  & -\Matrix{I}\end{pmatrix}
\begin{pmatrix} \MatrixSymbol{\alpha}\\ \Matrix{\W}_\rmt \end{pmatrix}\geq
\begin{pmatrix} \MatrixSymbol{\xi}_1^{} \\ \MatrixSymbol{\xi}_2\end{pmatrix},\hspace{0.1\textwidth}
\begin{pmatrix} \Matrix{I}  & \Matrix{0} \\ \Matrix{I} & -\Matrix{I}\end{pmatrix}
\begin{pmatrix} \MatrixSymbol{\beta}\\ \Matrix{\W}_\rmn \end{pmatrix}\geq\begin{pmatrix}
\Matrix{0} \\ \MatrixSymbol{\xi}_4\end{pmatrix},
\end{equation}
with the identity matrix $\Matrix{I}$ and $ \MatrixSymbol{\xi}_i$ corresponding to the right-hand sides in~\eqref{Eq_DiscreteRestrictions}.

Both inequalities are defined by full rank matrices.
Thus,  the following relations hold:
\begin{subequations}\label{Eq_SubstAB}
\begin{equation}\label{Eq_SubstA}
\begin{pmatrix} \MatrixSymbol{\alpha}\\ \Matrix{\W}_\rmt \end{pmatrix}=
\frac12\begin{pmatrix} \Matrix{I}  & \Matrix{I} \\ \Matrix{I}  & -\Matrix{I}\end{pmatrix}
\begin{pmatrix} \Matrix{y}_1 \\ \Matrix{y}_2\end{pmatrix},
\text{ with }\begin{pmatrix} \Matrix{y}_1 \\ \Matrix{y}_2\end{pmatrix}\geq
\begin{pmatrix} \MatrixSymbol{\xi}_1 \\ \MatrixSymbol{\xi}_2\end{pmatrix}
\end{equation}
and
\begin{equation}\label{Eq_SubstB}
\begin{pmatrix} \MatrixSymbol{\beta}\\ \Matrix{\W}_\rmn\end{pmatrix}=
\begin{pmatrix} \Matrix{I}  & \Matrix{0} \\ \Matrix{I} & -\Matrix{I}\end{pmatrix}
\begin{pmatrix} \Matrix{y}_3 \\ \Matrix{y}_4\end{pmatrix},
\text{ with } \begin{pmatrix} \Matrix{y}_3 \\ \Matrix{y}_4\end{pmatrix}\geq
\begin{pmatrix} \Matrix{0} \\ \MatrixSymbol{\xi}_4\end{pmatrix}.
\end{equation}
\end{subequations}
Thus, there is the same number of bound constraints as provided by the more general restrictions~\eqref{Eq_DiscreteRestrictions}.

Then, the discretized functional~\eqref{Eq_DiscreteEnergyV}  can be expressed in a general matrix form as
\begin{equation}\label{Eq_DiscreteGeneral}
\wt{\mathcal{K}}^k_h(\Matrix{y})
=\frac12\Matrix{y}^\top\Matrix{A}^k_h\Matrix{y}
-(\Matrix{b}^k_h)^\top\Matrix{y}+\Matrix{c}^k_h,
\qquad
\Matrix{y}\ge\MatrixSymbol{\xi}.
\end{equation}
The bound $\MatrixSymbol{\xi}$ is in fact determined by the constraints applied on $\Matrix{y}_j$, $j{=}1,2,3,4$ in~\eqref{Eq_SubstAB}.
The constrained minimum is denoted by $\Matrix{y}^k$.

The problem with standardly applied algorithms is that  the matrix \Matrix{A} might not necessarily be calculated in an explicit way.
The terms which arise from the  integrals in the first brackets in the right-hand side of~\eqref{Eq_DiscreteEnergyV} provide the energy and calculating the derivative with respect to the unknown \VectorMatrix{v} they provide a projected traction $\Matrix{M}\VectorMatrix{\P}$ with \Matrix{M} defined as in~\eqref{MatrixMeta}.
The projected traction can naturally be calculated from the SGBEM algorithm represented by the product $\Matrix{A}\Matrix{y}$   in equation~\eqref{Eq_DiscreteGeneral}.
Thus, each time the algorithm for minimization requires  a matrix-by-vector product, a linear system of  SGBEM is solved.
The influence matrices of SGBEM, however, are  calculated only once during the whole  solution process, as they are the same for all the iterations and time steps, considering only small displacements.

\begin{remark}[Solving QP problems]
\upshape
There are various ef\-fi\-cient numerical algorithms for solving QP, see \cite{dostal09B1} for a survey.
Here, we have used Conjugate Gradient (CG) based algorithms with bound constraints \cite{dostal09B1}.
The CG algorithm is not described explicitly  as it can be found together with all the necessary details of the constrained minimization, e.g., in reference~\cite{dostal09B1} and a scheme of the numerical implementation for a similar kind of minimization can also be found in~\cite{vodicka14A1}.
Additionally for larger problems preconditioning of CG methods is also necessary~\cite{dostal09B1} which together with a scalable algorithm can provide a robust solved for large problems~\cite{dostal10A1,dostal10A2}.
\end{remark}

\begin{remark}[Second-order cone programming]\label{rem-d=3}
\upshape
If $d\ge3$, $\mathcal{R}_1(\Z;\cdot)$ does not have a polyhedral graph as the tangent space is not one-dimensional and the friction can be anisotropic~\cite{PaThoSt09AFDS,Zmitrowicz06} and QP cannot be used.
Anyhow, such constraints are cones described by quadratic functions, and fall into the class of the \emph{second-order cone programming} (SOCP) for which efficient algorithms are still at disposal, e.g.~\cite{AliGol03SOCP,Stur02IIPM}.
\end{remark}

\section{Numerical examples}\label{Sec_NumExample}
The numerical procedure for the solution of friction contact problems described above is tested solving various sample problems.
The BE  code is based on the approach for treating   multi-domain problems derived in~\cite{vodicka07A1} which enables independent discretizations of solids on both sides of the contact zone $\GC$.
The solution of the contact problem is left upon a QP code, which uses standard algorithm of Modified Proportioning with Reduced Gradient Projections (MPRGP) based on the schemes in~\cite{dostal09B1}.

\subsection{Description of the contact layouts}\label{Sec_NumExampleDescription}

In the examples presented, friction contact problems  including typical receding and conforming contacts with proportional and non-proportional loading  are solved.

\begin{figure}[ht]
\centering
\PictureFile{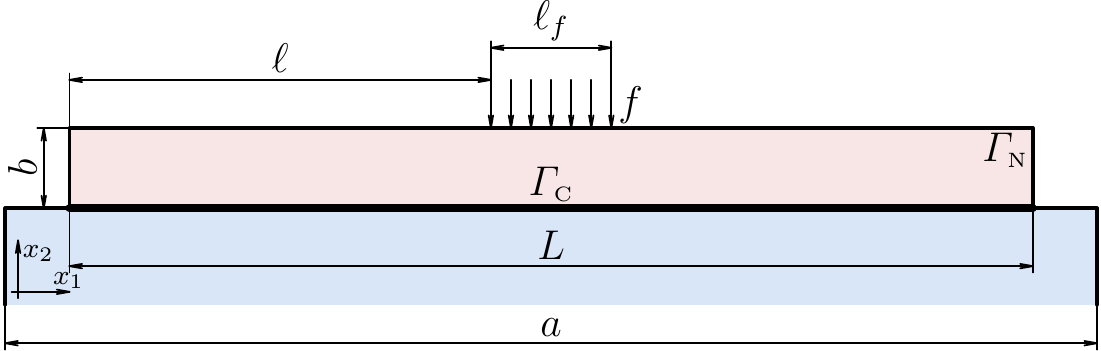}
\caption{Receding contact layout.}\label{Fig_Receding}
\end{figure}

In the first example, the geometry and load arrangement correspond to a typical receding contact problem.
The plane strain problem layout is shown in Figure~\ref{Fig_Receding}.
The bottom block is not plotted entirely as it has, in fact, a square cross-section.
The dimensions used in the picture are: $L{=}160\,$mm, $b{=}10\,$mm, $a{=}200\,$mm, $\ell{=}75\,$mm, $\ell_{\Vector{f}}{=}10\,$mm.
The visco-elastic properties of the blocks are $E{=}4{\times}10^3$MPa, $\nu{=}0.35$, and the relaxation-time parameter  is $\chi{=}10^{-3}$s.

The contact is characterized by the friction coefficient and the normal compliance: the Coulomb friction coefficient is $\mu{=}0.8$, the normal stiffness is $k_g{=}4{\times}10^5$MPa\,mm$^{-1}$.

In this example, three  BE meshes are used, each with a particular time-step.
The coarsest spatial BE mesh is refined at the contact zone $\GC$,  at the short sides of the layer and in the area of the applied load,  whereas it is  more or less uniform, with the typical element length $\ell_\text{e}{=}20\,$mm, along the remaining part of the boundary.
The smallest element length is   $\ell_\text{e\,min}{=}4\,$mm in the central part of   $\GC$ for $x_1{\in}[80;120]$.
The time-step is $\TAU{=}1{\times}10^{-3}$s in this coarsest mesh.
This discretization is denoted $N{=}10$ according to the number of elements in the uniform part of the contact zone mesh.
The refined meshes are denoted subsequently as $N{=}20$ and $N{=}40$: the lengths of all the boundary elements   and also the time steps  are divided by two with respect to the previous coarser discretization.

The vertical load \Vector{f} is applied incrementally in time-steps such that the maximum applied pressure is $f_\text{max}{=}0.5$MPa.
The bottom face of the bottom square block is fixed.

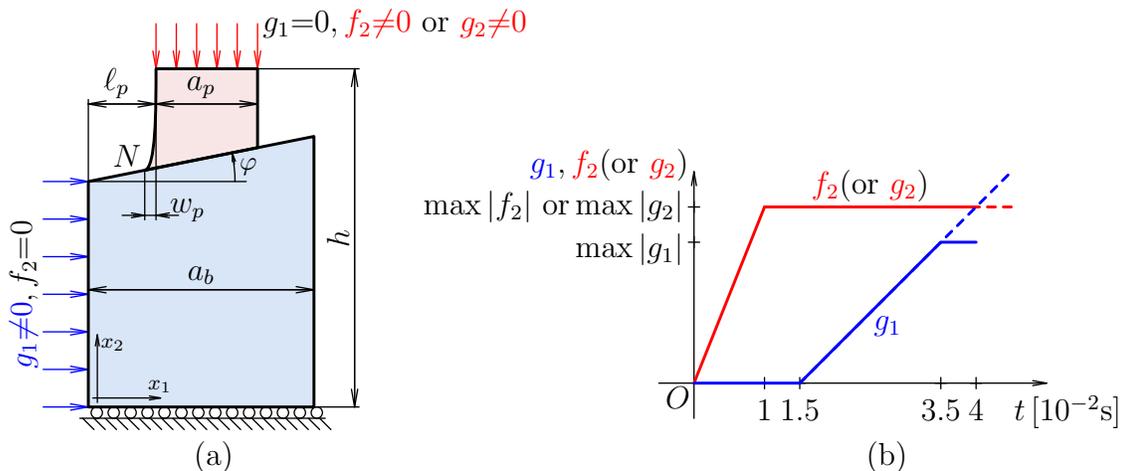
\begin{figure}[ht]
\centering
\settoheight{\unitlength}{\PictureFile{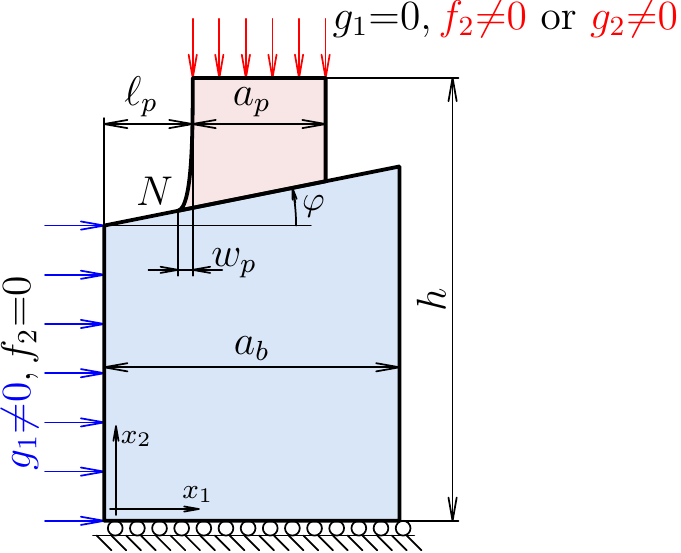}}
\begin{picture}(0,1)
\put(0,0.025){\makebox(0,0)[br]{\PictureFile{Geometry03}}}
\put(-0.25,0.025){\makebox(0,0)[bl]{\PictureFile{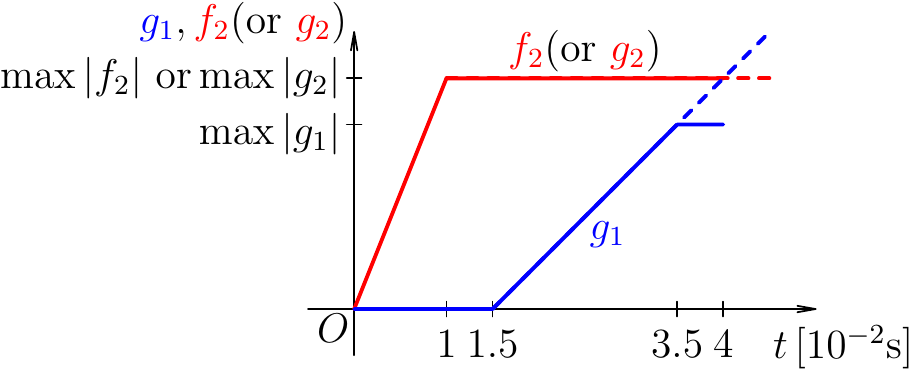}}}
\put(-0.7,0){\makebox(0,0)[tr]{(a)}}
\put(0.8,0){\makebox(0,0)[tl]{(b)}}
\end{picture}
\caption{
The layout for the second and the third examples: (a) geometry and (b) time dependence of external loading.
}\label{Fig_Block}
\end{figure}

The layout of the plane strain problem in the second and the third examples is shown in Figure~\ref{Fig_Block}.
Geometrically, the problem is a conforming contact problem if $\varphi{=}0$ and $g_1{=}0$.
Nevertheless, the prescribed load enables the lower block to slide if the frictional forces are exceeded and then the upper block may detach partially or totally obtaining again receding contact.

The load is applied non-proportionally according to the scheme shown in Figure~\ref{Fig_Block}.
First, the vertical load $f_2$ or $g_2$ is  increasing in time-steps up to its maximum value, then the structure relaxes for a while keeping the load constant.
Notice that, the upper side of the punch is prevented from horizontal moving in any case.
Subsequently, the horizontal pushing $g_1$ starts increasing up to the maximum prescribed displacement.
The tangential forces along the vertical face loaded by $g_1$ are kept zero so that prior to initializing the horizontal load both bottom and left face of the bottom block are constrained equally by simple support conditions.

For the second example, the dimensions used in  Figure~\ref{Fig_Block} are: $a_b{=}200\,$mm, $a_p{=}100\,$mm, $\ell_p{=}50\,$mm, $h{=}300\,$mm, $\varphi{=}0$, $w_p{=}0\,$mm.
The physical parameters are the same as in the previous example, with exceptions of Young modulus of the bottom block which is $E_b{=}4E$ and the friction coefficient which is $\mu{=}0.2$.

The used BE mesh  is uniform, with element length $\ell_\text{e}{=}5\,$mm, far from $\GC$,  and is refined close to the end points of $\GC$ due to expected stress singularity at these points, where the smallest element length is   $\ell_\text{e\,min}{=}0.11\,$mm.
There are 80 elements on each side  of the contact zone $\GC$ (i.e. in each solid $A$ and $B$).
The time-step $\TAU{=}2.5{\times}10^{-4}$s.

The vertical load is $f_2$ and its maximum value is $f_\text{max}{=}1$MPa.
Similarly, the maximum prescribed horizontal displacement $g_1$ is  $g_\text{max}{=}0.1\,$mm, which is still less then the minimum element length $\ell_\text{e\,min}$.

For the third example, the modified geometrical data are $a_p{=}25\,$mm, $\ell_p{=}80\,$mm, $\varphi{=}\arctan\frac12$, $w_p{=}5\,$mm.
A `nose' of the punch (the point $N$ in Figure~\ref{Fig_Block}) is introduced in order to eliminate stress singularity at this point, see~\cite{comninou76}, which could lead to incision of  the punch into the lower block therein and prevent it even from partial sliding (slipping).
This singularity may be rather strong for higher friction coefficients when trying to slide  the bottom block  to the right, see~\cite{comninou76} for details.
The angle of the punch domain at the nose is about $5^{\circ}$.
The physical parameters are still kept the same with an exception of the friction coefficient $\mu$ and the relaxation-time parameter $\chi$.
Both are varied to affect the behaviour of the system: $\mu$ varies between $0.2$ and $1.1$, $\chi$ may change from $0$ to $0.01$s.

The main features of the BE mesh used in the third example are: element length $\ell_\text{e}{=}20\,$mm far from $\GC$,  and refined close to the end points of $\GC$ due to possible stress concentrations at these points, where the smallest element length is   $\ell_\text{e\,min}{=}0.25\,$mm.
There are 46 elements on each side of the contact zone $\GC$ whose lengths close to the endpoints are  $\ell_\text{e\,min}$ and in the central part of $\GC$ the maximal length is four times more.
The initial time-step $\TAU{=}1{\times}10^{-3}$s.

The load is also applied in slightly different way, though following the same time scheme as shown in~\ref{Fig_Block}(b).
The vertical load now is $g_2$ and its maximum value is $g_\text{max}{=}0.02\,$mm.
Similarly, the maximum prescribed horizontal displacement $g_1$ causing in all the cases detachment of the blocks is  $g_\text{max}{=}0.21\,$mm, which is still less then the minimum element length $\ell_\text{e\,min}$.

In this example, an adaptive algorithm based on the discrete energy inequality~\eqref{Eq_EbilanceDiscr} is verified.
A parameter $\epsilon$ is set and if the difference between the left-hand and right-hand sides is greater than $\epsilon$, the time step is multiplied by a half, if the
difference is less than $0.1\epsilon$, the time step is multiplied by two, otherwise it is kept the same.
Simultaneously, a minimal time-step length is set to $\TAU_{\min}{=}10^{-6}$s.
In calculation, $\epsilon$ is chosen in the interval form $64\,\upmu$J to
$1\,\upmu$J.

The way of loading is the same as in the second example only the horizontal loading does not stop until the two bodies separate (the dashed line in the right drawing of Figure~\ref{Fig_Block}).

\subsection{Results for the receding contact}\label{Sec_ResultsReceding}

\begin{figure}[ht]
\setlength{\unitlength}{1\textwidth}
\centering
\begin{picture}(0.93,0.32)
\put(0,0){\makebox(0,0)[bl]{\PictureFile{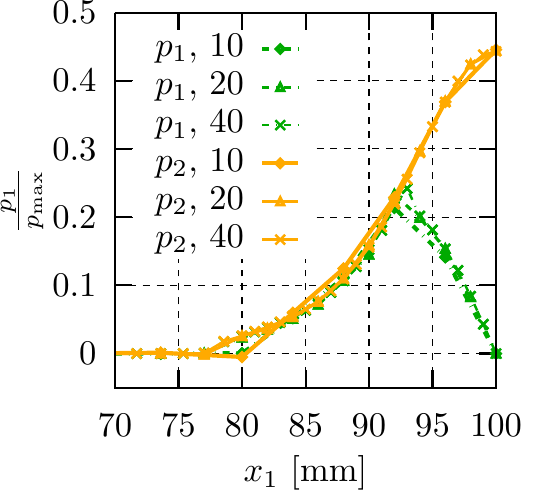}}}
\put(0.31,0){\makebox(0,0)[bl]{\PictureFile{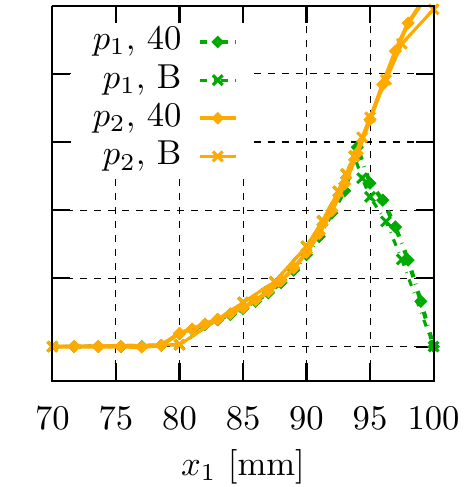}}}
\put(0.58,0){\makebox(0,0)[bl]{\PictureFile{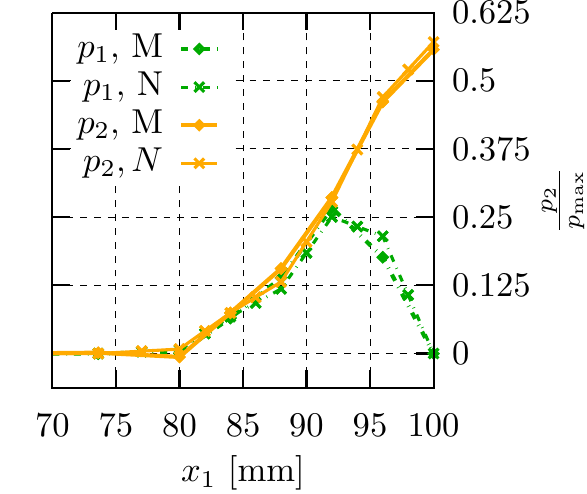}}}
\put(0.09,0.32){\makebox(0,0)[t]{(a)}}
\put(0.36,0.32){\makebox(0,0)[t]{(b)}}
\put(0.63,0.32){\makebox(0,0)[t]{(c)}}
\end{picture}
\caption{Distributions of contact stresses ($x_1{=}100\,{\rm mm}$ is the centre of the specimen in Fig.\,\ref{Fig_Receding}): (a) the three discretizations with $N$ introduced in Section~\ref{Sec_NumExampleDescription}; (b) comparison of the finest discretisation $N{=}40$ with the approach in~\cite{blazquez98} (B); (c) matching (M) and non-matching (N) discretizations.}\label{Fig_StressReceding}
\end{figure}
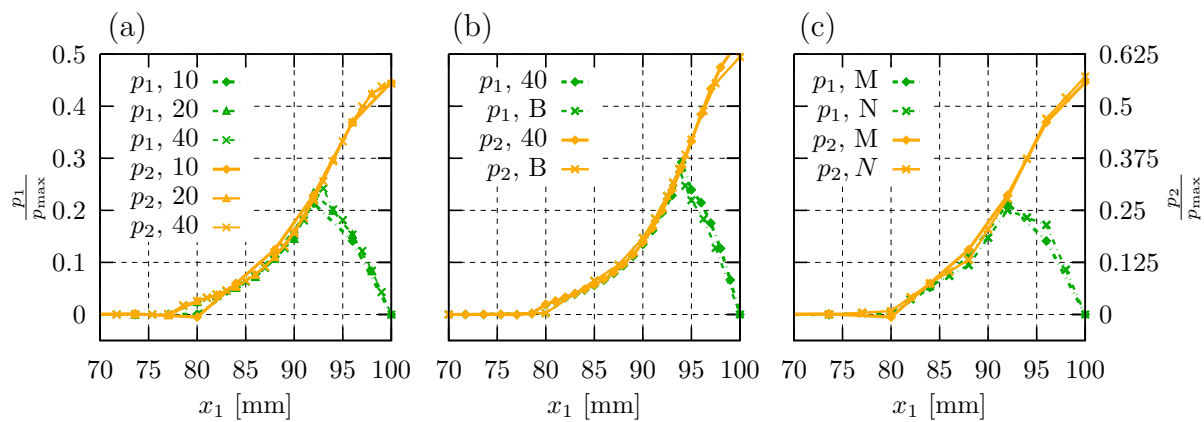

In this example we demonstrate  how the solution of the contact problem converges.
The graphs in Figure~\ref{Fig_StressReceding} show only a part of the left half of the interface due to symmetry.
As $x_1{=}0\,$mm pertains to the left edge of the bottom block, the value of $x_1{=}100\,$mm is exactly the midpoint of the contact zone.
Generally notice that scaling of the axes by the coefficient of friction  shows the normal and tangential component to be equal in the slip zone.


Figure~\ref{Fig_StressReceding}(a) shows the calculated contact tractions obtained  by the three defined meshes after the whole load has been applied.
The numerical results approximate the actual solution very well which can also be seen in Figure~\ref{Fig_StressReceding}(b).
Here, the present results (for the finest mesh $N{=}40$) are compared to those obtained by a classical contact algorithm of Bl\'azquez et al.~\cite{blazquez98}.
An excellent agreement between the results of both approaches is achieved, though different contact algorithms and quite coarse meshes employed cause small differences.
It should be mentioned that instead of  the present distributed loading, a point load in the middle of the upper layer was used in the calculation to be conformable with in~\cite{blazquez98}. 
Finally in Figure~\ref{Fig_StressReceding}(c), a  comparison of the stress solutions for matching and a non-matching discretizations in the contact zone is carried out, enabled by the computer code based on~\cite{vodicka07A1}..
To this end, the coarsest mesh from Figure~\ref{Fig_StressReceding}(a) is modified splitting each element  into two equal elements along one side of the contact zone.
A nice agreement between the matching mesh and the finer side of the non-matching mesh can be observed, with only small differences between both results, acceptable  for the used BE meshes.

\begin{figure}[ht]
\centering
\PictureFile{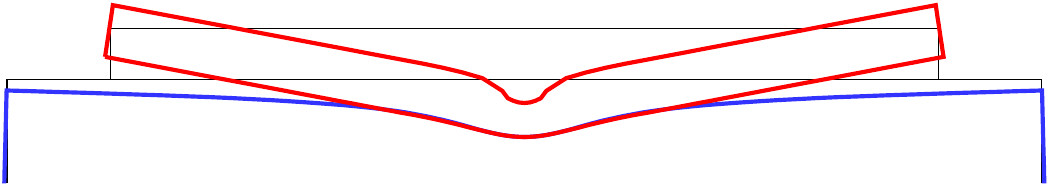}
\caption{Deformed (magnified ${\times}5000$) and undeformed configuration of the layer with receding contact from Fig.\,\ref{Fig_Receding} at the end of the loading process.}\label{Fig_DefReceding}
\end{figure}

Finally, the deformed configuration at the maximum load is plotted  in Figure~\ref{Fig_DefReceding} to demonstrate the receding character of the problem and also that the penalization of the normal contact condition does not cause any significant interpenetrations of the contacted bodies.

\subsection{Results for the conforming contact}\label{Sec_ResultsBlock}

\begin{figure}[ht]
\centering
\setlength{\unitlength}{1\textwidth}
\begin{picture}(0.9,0.35)
\put(0,0.25){\makebox(0,0)[bl]{\PictureFile{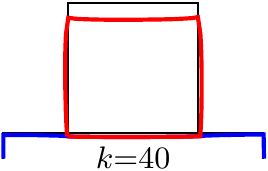}}}
\put(0,0.14){\makebox(0,0)[bl]{\PictureFile{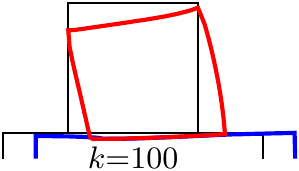}}}
\put(0,0.03){\makebox(0,0)[bl]{\PictureFile{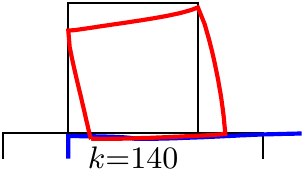}}}
\put(0.22,0.025){\makebox(0,0)[bl]{\PictureFile{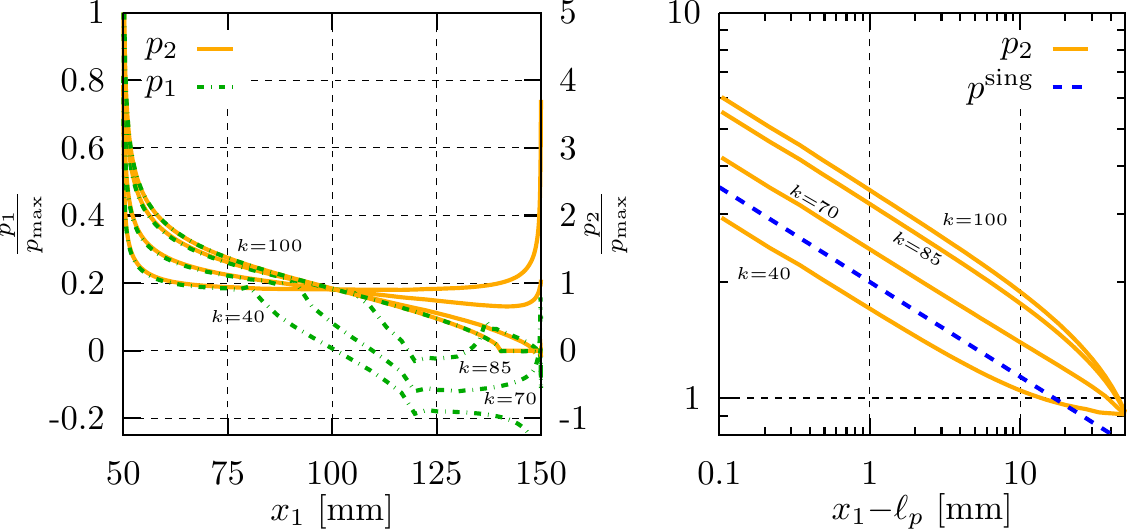}}}
\put(0.07,-0.01){\makebox(0,0)[bl]{(a)}}
\put(0.4,-0.01){\makebox(0,0)[bl]{(b)}}
\put(0.77,-0.01){\makebox(0,0)[bl]{(c)}}
\end{picture}
\caption{Conforming contact,  Fig.~\ref{Fig_Block} with $w_p{=}0$ and $\varphi{=}0$: (a)  deformed configurations (magnified ${\times}500$) at three selected instants of the loading process; (b) evolution of the spatial distribution of contact stresses, $x_1{=}100\,$mm lies in the midpoint of the contact zone; (c) a detail at the stress singularity point $x_1{=}\ell_p$.}\label{Fig_ResBlock}
\end{figure}
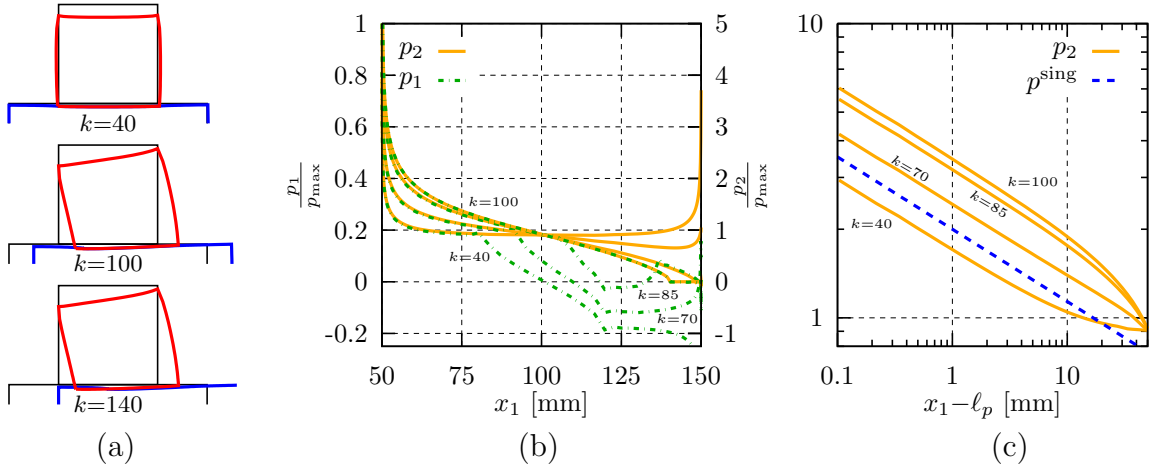

In this example we  test the numerical procedure for  a non-proportional loading with results in  Figure~\ref{Fig_ResBlock}.
First, Figure~\ref{Fig_ResBlock}(a) shows  the deformation of the punch and top part of the bottom block at three time-steps:
$k{=}40$, i.e.\ $t{=}0.01$s, the vertical load in Figure~\ref{Fig_Block} stops  its increase, no lateral loading;
$k{=}100$, i.e.\ $t{=}0.025$s, the frictional forces are exceeded by the horizontal load, the bottom block strarts to move;
$k{=}140$, i.e.\ $t{=}0.04$s${=}T$, the end of loading, the bottom block has moved with respect to the punch.
The shown magnified displacements are in fact smaller than the minimum element  length so that the displayed interpenetration of the bodies is only fictitious.

The contact stress behaviour is  demonsrated in Figure~\ref{Fig_ResBlock}(b), where the curves correspond to the aforementioned load steps $k{=}40$ and $k{=}100$ and the intermediate time-steps $k{=}75$ and $k{=}80$ to better understand a quite complicated evolution: from typical conforming contact with two stress singularities to a sliding  block with only one singularity, here  at the left corner of the punch.
The changes between stick in the central part to partial sliding in the outer part of the contact zone cause some non-smooth changes in the distribution of the tangential tractions, and partially can be eliminated by refining the mesh because these non-smooth changes occur when the actual transition point between the stick and the slip zones  lies in the middle of an element of the given discretization.

Finally, the graph plotted in log-log scale in Figure~\ref{Fig_ResBlock}(c) documents that the stresses near the singularity  point follow the power law  $p^{\text{sing}}{=}K(x_1{-}\ell_p)^{-0.24582}$, where the singularity exponent, depending on the coefficient of friction not on the applied load, has been evaluated according to~\cite{comninou76}.

\subsection{Results for the skewed punch}\label{Sec_ResultsSkew}

In this example we  test the role of the coefficient of friction  on the quality of the solution of the contact problem, the role of viscosity through the time relaxation parameter $\chi$ and the influence of time adaptivity by setting prescribed difference $\epsilon$ in the inequality~\eqref{Eq_EbilanceDiscr}.

The first set of pictures presents  the energy residuum $\Delta E$ obtained as  the difference between the right- and the left-hand sides of~\eqref{Eq_EbilanceDiscr}.
It also shows the evolution of the  total reaction forces at a part of the boundary: $F_1$ acts along the left edge of the block, the vertical forces there are zero; $F_2$ act along the bottom
edge of the block, the horizontal forces vanish at this edge.
In each of the three tests one parameter was changing to obtain a few graphs and to asses the influence of each particular parameter.

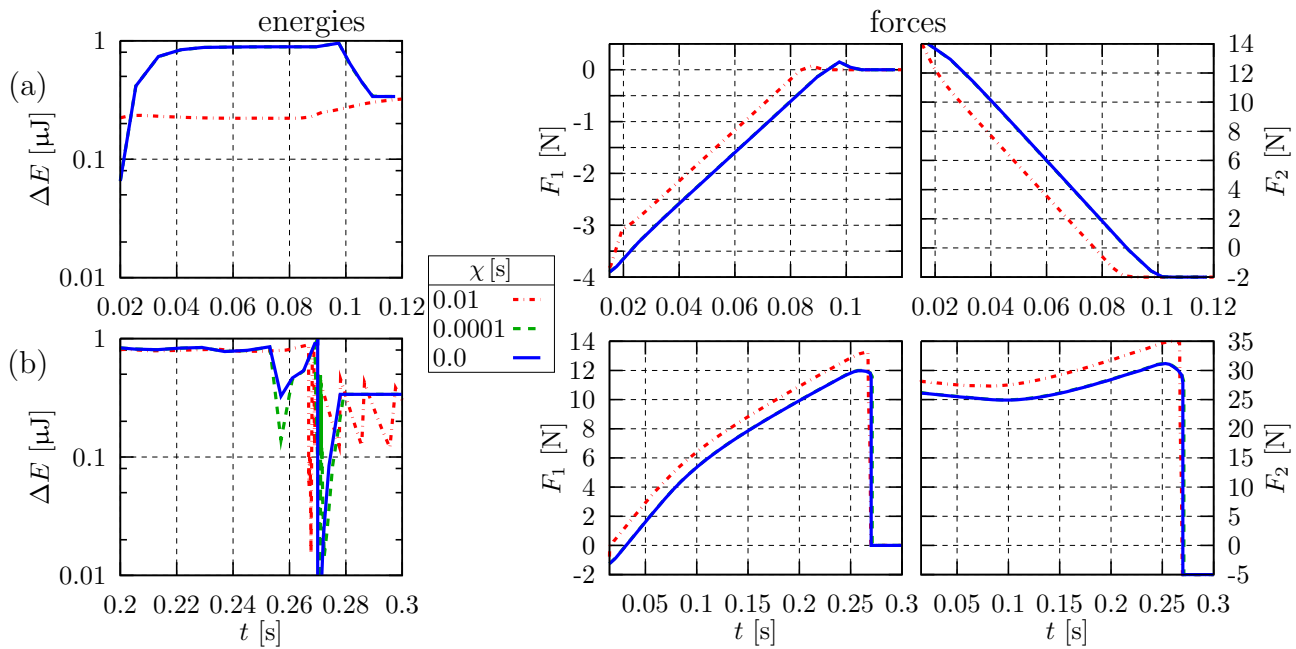
\begin{figure}[ht]
\centering
\settoheight{\unitlength}{\PictureFile{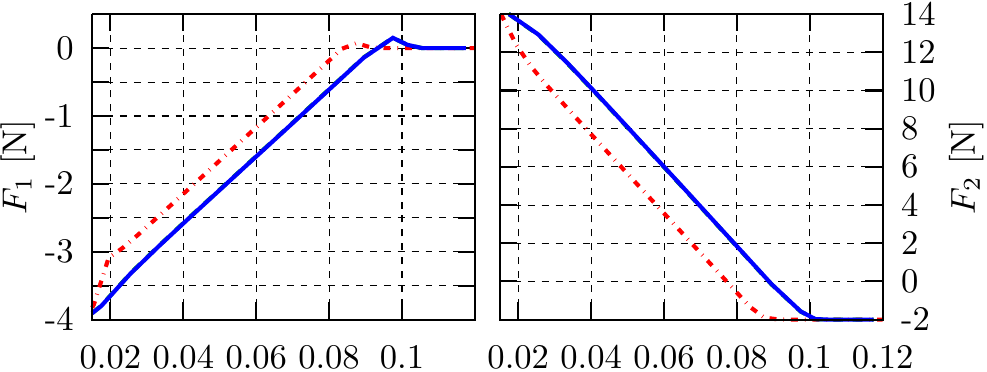}}
\begin{picture}(4.5,2.25)
\put(0,0.0){\makebox(0,0)[bl]{\PictureFile{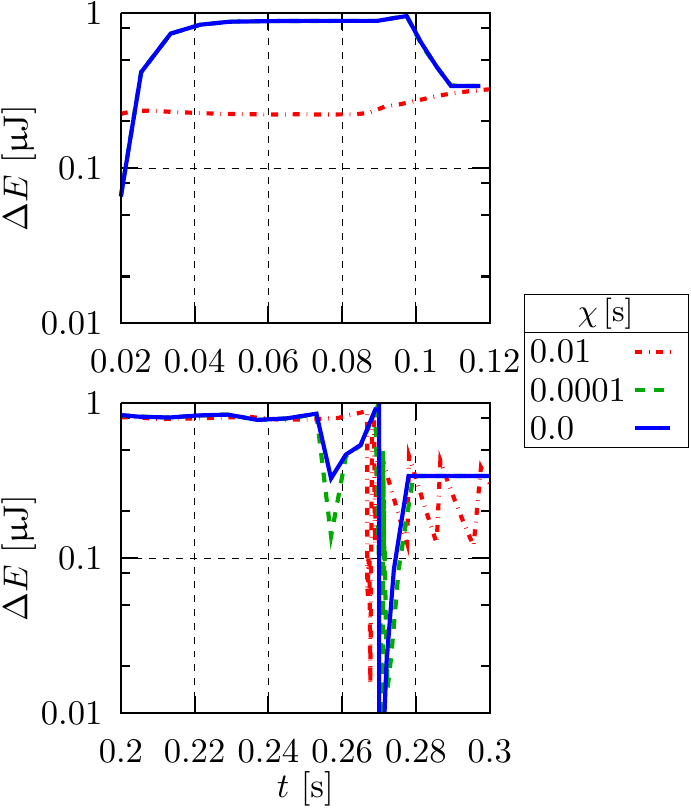}}}
\put(1.8,1.17){\makebox(0,0)[bl]{\PictureFile{Skew_Forces02}}}
\put(1.8,0.0){\makebox(0,0)[bl]{\PictureFile{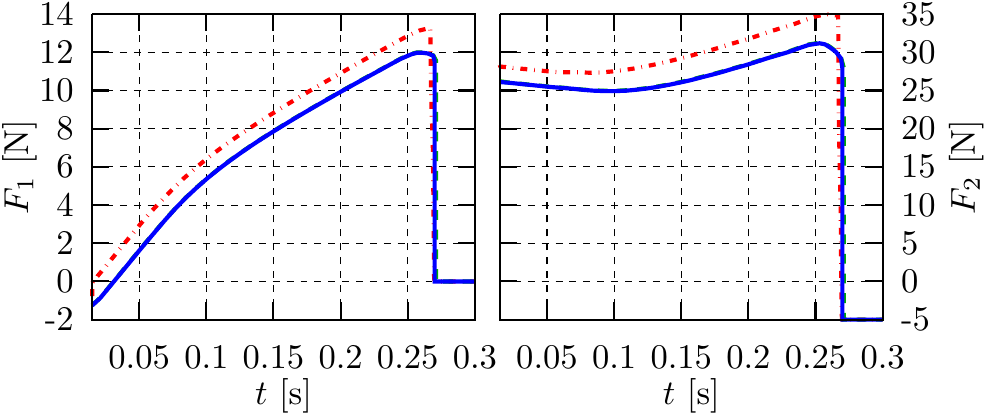}}}
\put(1,2.25){\makebox(0,0)[tc]{energies}}
\put(3.1,2.25){\makebox(0,0)[tc]{forces}}
\put(0,1.98){\makebox(0,0)[cc]{(a)}}
\put(0,1){\makebox(0,0)[cc]{(b)}}
\end{picture}
\caption{Time evolution of the energy residuum and total forces at the constraint boundaries of the bottom block for various values of time-relaxation parameters $\chi$ and for the residual tolerance $\epsilon{=}1\,\upmu${\rm J}: (a) small friction $\mu{=}0.2$, (b) large friction $\mu{=}1.1$. 
The response for the small viscosity $\chi{=}0.0001\,${\rm s} essentially coincides with the inviscid case $\chi={0}\,${\rm s} (for which no theoretical supportin Appendix~\ref{Sec_Theory} is given, however.
Notice also different time ranges used in the cases (a) and (b).
}\label{Fig_SkewEFTR}
\end{figure}

Figure~\ref{Fig_SkewEFTR} compares evolutions of $\Delta E$ and total forces for different values of $\chi$. 
The residual tolerance is the finest used, i.e.\ $\epsilon{=}1\,\upmu$J, and this upper bound is set as the maximum of the vertical axis range.
Figures~\ref{Fig_SkewEFTR}(a) and (b) are obtained for the  smallest and the greatest values  of the coefficient of friction, $\mu{=}0.2$ and $\mu{=}1.1$, respectively. 
There is no significant difference between both friction cases when observing the convergence with respect to $\chi$, besides the energy residuum for larger friction where a sharp peak appears which corresponds to the time instant when the punch separates from the foundation.
Nevertheless, for small viscosities and for switched-off viscosity the graphs for energy resuduum and total forces, too, are similar and close to each other (in fact  they can be distinguished only energy plot in  Figure~\ref{Fig_SkewEFTR}(b)).
The graphs show only a part of loading history after the horizontal push has been initiated.
There would be some differences in total forces  due to different viscosity parameters, if the time instants of load changes (see Figure~\ref{Fig_Block}) were presented.
Of course, the overall behaviour of the total force is different, as in the small-friction case the punch slides on the foundation in the presented part of the load history and in the large-friction case there is a stick zone in the contact so that the signs of the tangential tractions are opposite in the two cases.
Thus, during the stage of pushing the bottom block, the  traction decreases continuously until separation in the former case.
In the latter case, the tangential traction increases while the friction is capable to bear the loading, which is then followed by an abrupt separation of the bodies shown as a jump in the traction force distribution.
This seems to be a very interesting behaviour, which occurs for large friction only.

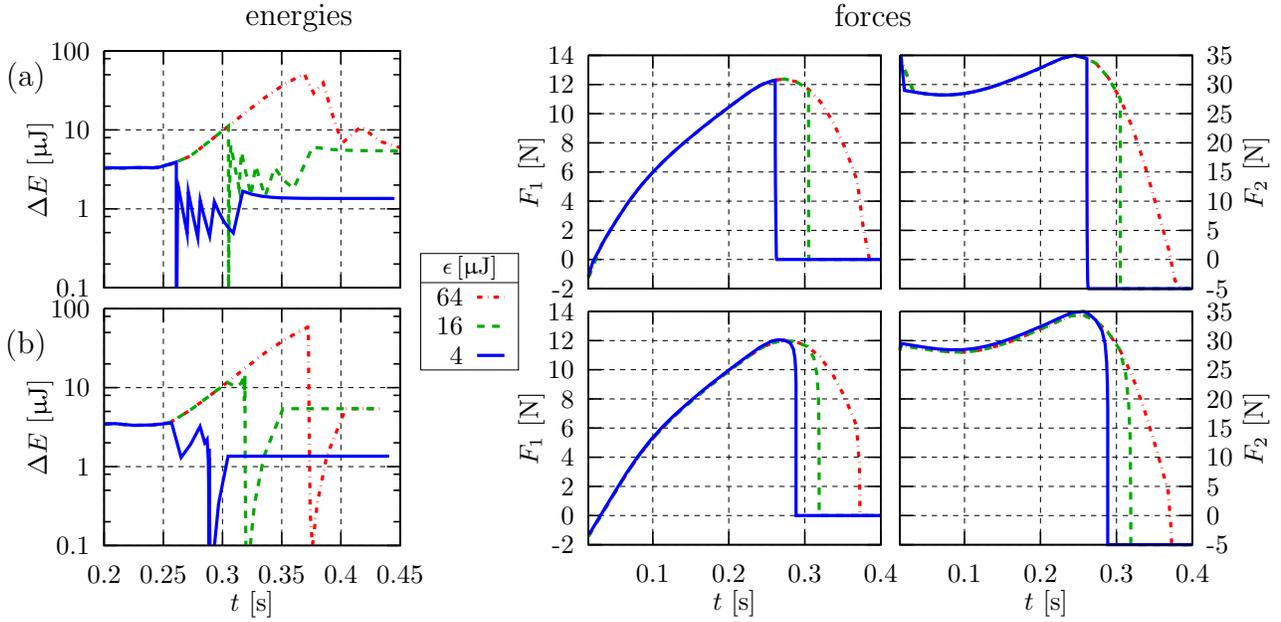
\begin{figure}[ht]
\centering
\settoheight{\unitlength}{\PictureFile{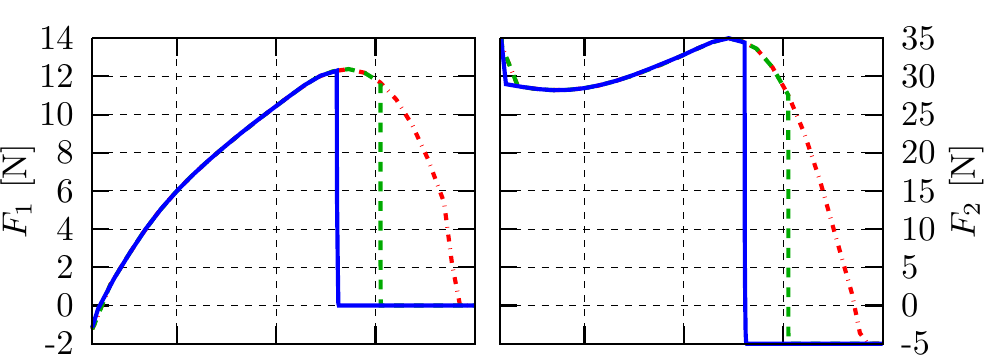}}
\begin{picture}(4.5,2.25)
\put(0,0.0){\makebox(0,0)[bl]{\PictureFile{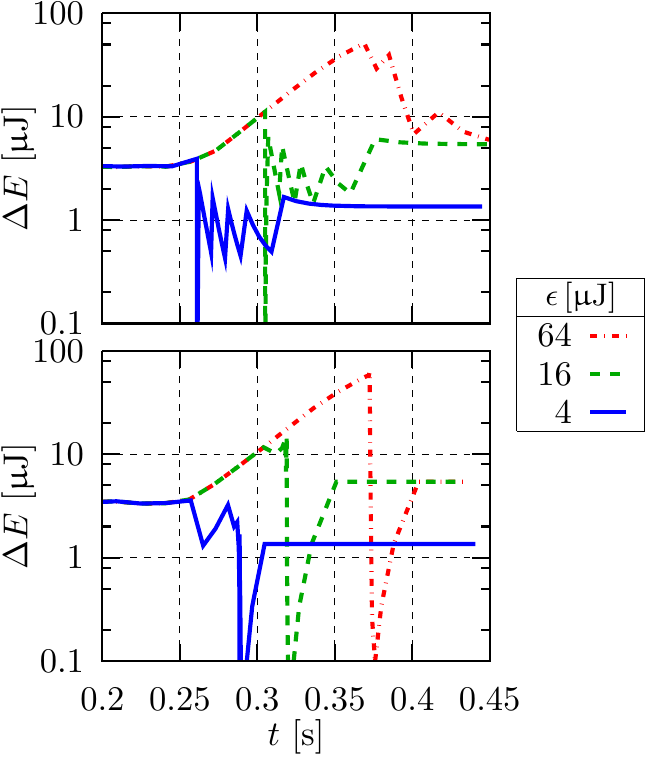}}}
\put(1.8,1.17){\makebox(0,0)[bl]{\PictureFile{Skew_Forces00}}}
\put(1.8,0.0){\makebox(0,0)[bl]{\PictureFile{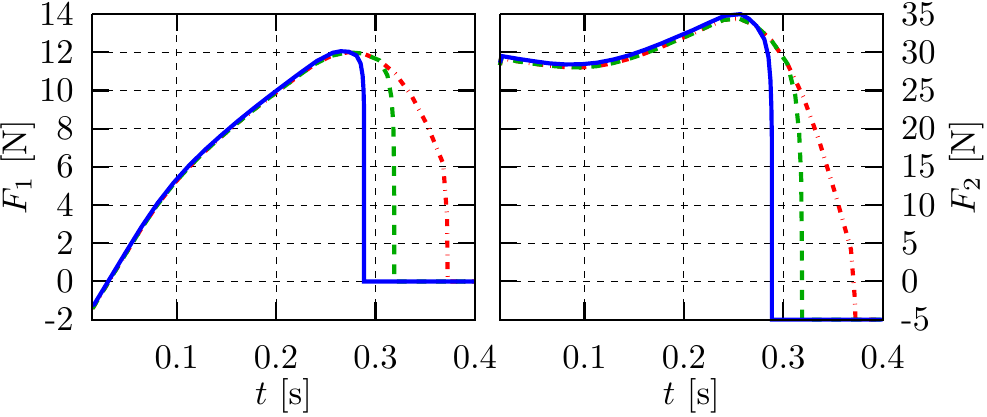}}}
\put(1,2.25){\makebox(0,0)[tc]{energies}}
\put(3.1,2.25){\makebox(0,0)[tc]{forces}}
\put(0,1.98){\makebox(0,0)[cc]{(a)}}
\put(0,1){\makebox(0,0)[cc]{(b)}}
\end{picture}
\caption{Time evolution of the energy residuum and total forces at the constraint boundaries of the bottom block for various tolerances $\epsilon$ and for the friction coefficient $\mu{=}1.1$:
(a) larger viscosity $\chi{=}0.01\,${\rm s}, (b) smaller viscosity $\chi{=}0.0001\,${\rm s}.}\label{Fig_SkewEFED}
\end{figure}

Therefore, when testing the influence of the parameter $\epsilon$, only graphs for the greatest coefficient of friction  $\mu{=}1.1$ are shown.
Figure~\ref{Fig_SkewEFED} compares evolutions of $\Delta E$ and total forces for various values of $\epsilon$.
Figures~\ref{Fig_SkewEFED}(a) and (b) show $\Delta E$ for the most viscous case considered with $\chi{=}0.01$s and the less viscous case with $\chi{=}0.0001$s, respectively.
Only the most interesting detail of the time instant of block separation is shown in the energy graphs.
In any case (energy or force), a convergence  of the block-separation instant for diminishing $\epsilon$ is evident and also an abrupt change from contact to separation can be clearly seen. 
The maxima of $\Delta E$ are perfectly bounded by the pertinent residual tolerances $\epsilon$.

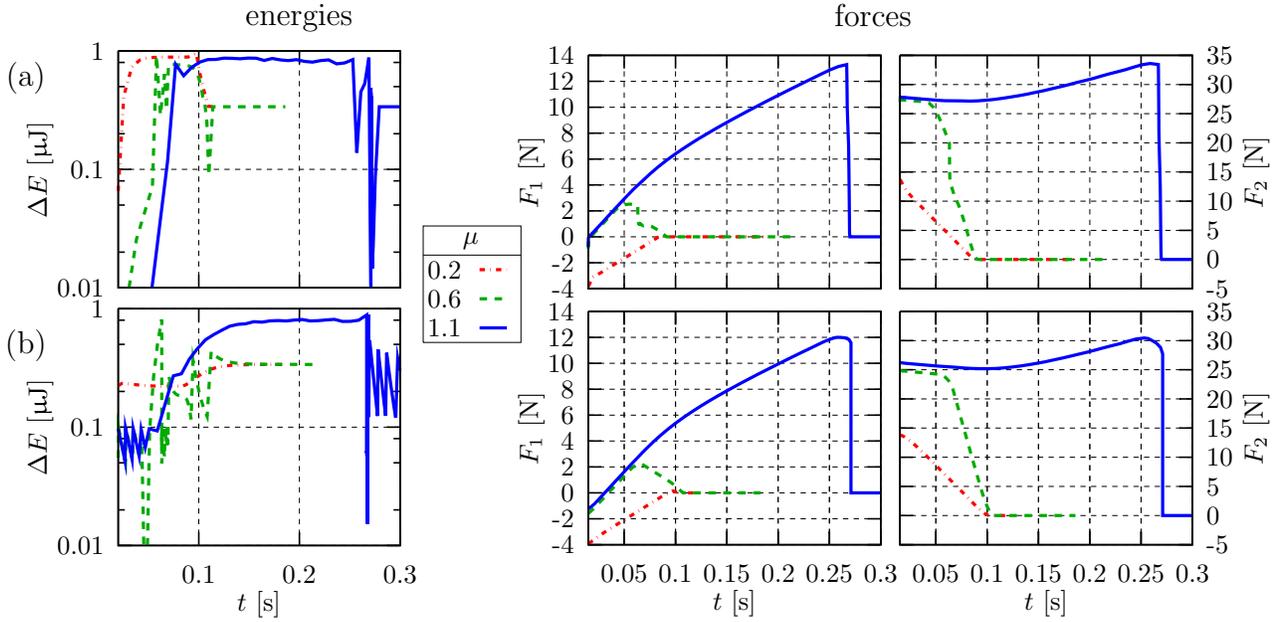
\begin{figure}[ht]
\centering
\settoheight{\unitlength}{\PictureFile{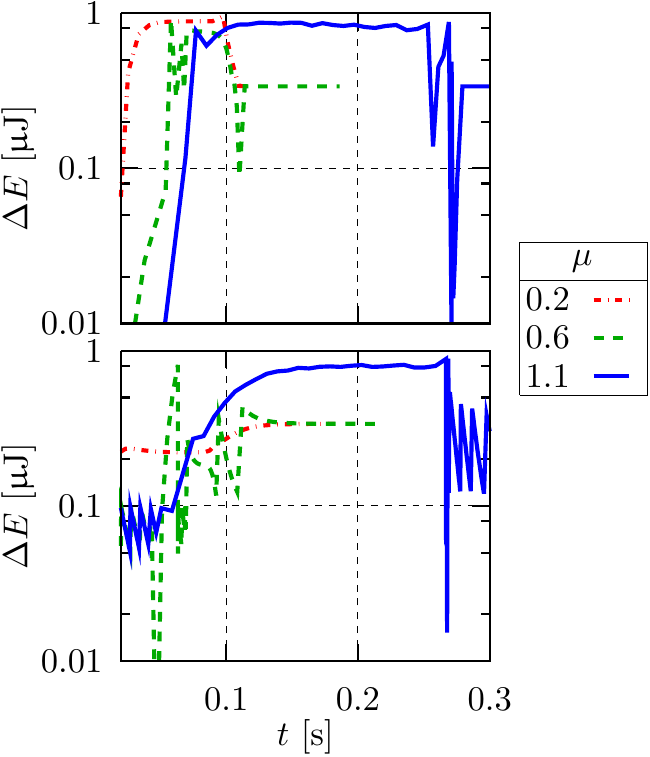}}
\settoheight{\unitlength}{\PictureFile{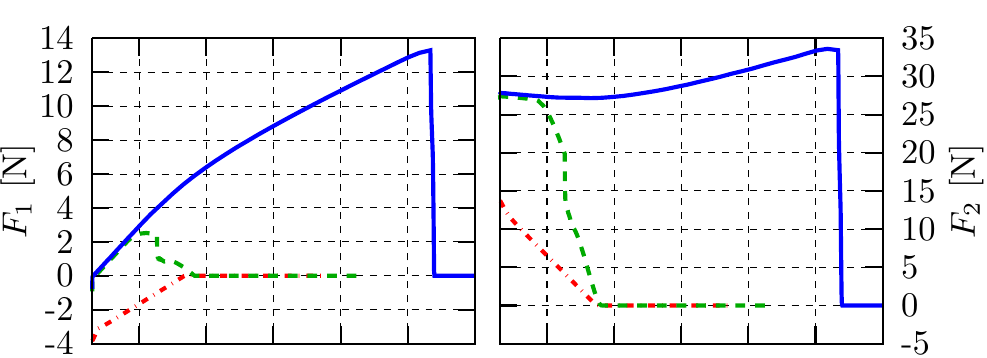}}
\begin{picture}(4.5,2.25)
\put(0,0.0){\makebox(0,0)[bl]{\PictureFile{Skew_Energies02}}}
\put(1.8,1.17){\makebox(0,0)[bl]{\PictureFile{Skew_Forces04}}}
\put(1.8,0.0){\makebox(0,0)[bl]{\PictureFile{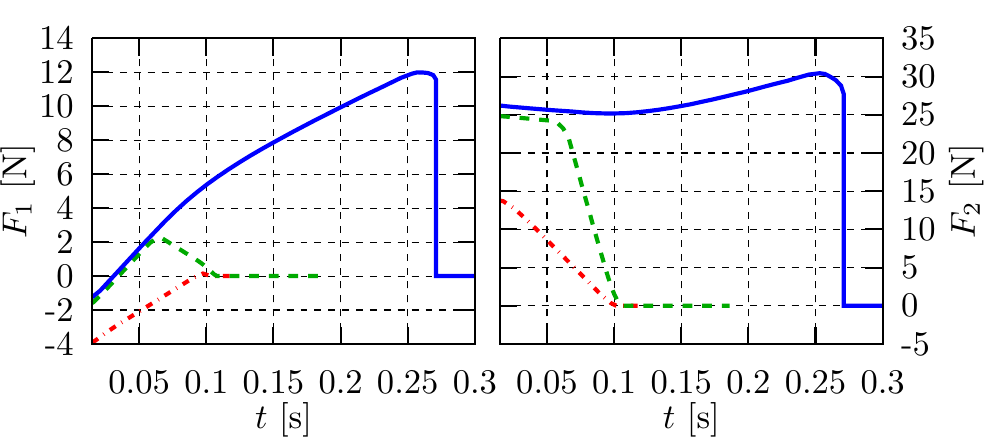}}}
\put(1,2.25){\makebox(0,0)[tc]{energies}}
\put(3.1,2.25){\makebox(0,0)[tc]{forces}}
\put(0,1.98){\makebox(0,0)[cc]{(a)}}
\put(0,1){\makebox(0,0)[cc]{(b)}}
\end{picture}
\caption{Time evolution of the energy residuum and total forces at the constraint boundaries of the bottom block for various friction coefficients $\mu$ and for the tolerance $\epsilon{=}1\,\upmu${\rm J}:
(a) larger viscosity $\chi{=}0.01\,${\rm s}, (b) smaller viscosity $\chi{=}0.0001\,${\rm s}.
}\label{Fig_SkewEFMU}
\end{figure}

Finally, Figure~\ref{Fig_SkewEFMU} compares $\Delta E$ and total forces  for different values of $\mu$.
The residual tolerance is again $\epsilon{=}1\,\upmu$J so that the maximum of the vertical axis range was set to this value.
Figures~\ref{Fig_SkewEFMU}(a) and (b) corresponde again to  $\chi{=}0.01$s and  $\chi{=}0.0001$s, respectively.
The increasing friction requires naturally more load, i.e.\ larger prescribed displacements applied on the left side of the block, to separate the two bodies. 
Thus, the jumps of the energy residuum appear later in time.
At small friction, there is a continuous decrease of the resultant force ending by separation of the two blocks.
The greatest friction provides a continuous increase of the resultant force terminated by an abrupt jump to zero at the moment of the block separation.
Somewhere in between, here represented by the value $\mu{=}0.6$, the forces show a mixed behaviour, including even small jumps for  the most viscous case.

\begin{figure}[ht]
\centering
\settoheight{\unitlength}{\PictureFile{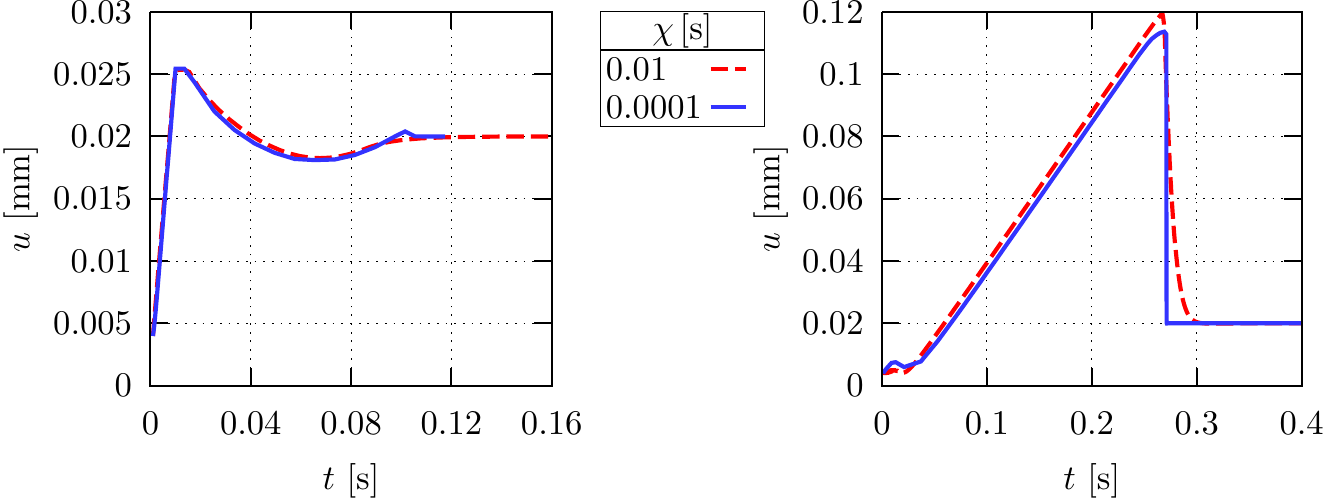}}
\begin{picture}(0,1)
\put(0,0.025){\makebox(0,0)[bc]{\PictureFile{Skew_Jump00}}}
\put(-0.2,0){\makebox(0,0)[tr]{(a) small friction $\mu{=}0.2$}}
\put(0.45,0){\makebox(0,0)[tl]{(b) large friction $\mu{=}1.1$}}
\end{picture}
\vspace*{.3em}
\caption{Time evolution of the total displacement of the nose $N$ from Fig.~\ref{Fig_Block}(a) for both viscosities $\chi=0.01\,${\rm s} and $0.0001$\,{\rm s}.
}\label{Fig_SkewJumpDef}
\end{figure}
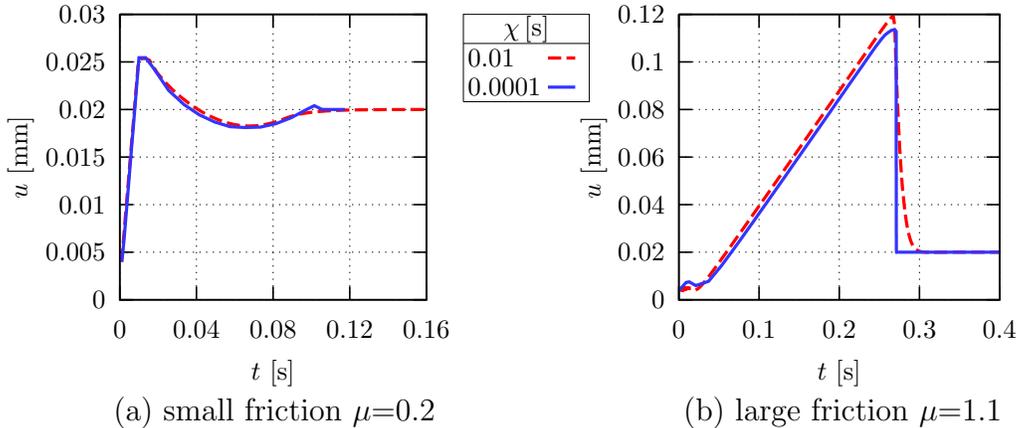

The jump behaviour  associated to large friction can be also seen in deformations. 
Figure~\ref{Fig_SkewJumpDef} shows evolution of the total displacement $u{=}\sqrt{u_1^2{+}u_2^2}$ at the point $N$ (of the upper block, see Figure~\ref{Fig_Block}(a)).
Considering various friction coefficients and viscosity parameters, we can conclude that, with large friction, the point jumps (i.e. there is a sudden change of the deformation state) when the frictional forces are not capable to keep the bodies in contact.
The body then returns to its initial configuration (the speed depends on the relaxation time $\chi$) shifted by the the prescribed displacement: $0.02\,$mm in vertical direction, see Section~\ref{Sec_NumExampleDescription}.

\begin{figure}[ht]
\centering
\settoheight{\unitlength}{\DeformationFile{02}{0001}{17}}
\begin{picture}(0,6)
\put(-0.01,4.05){\makebox(0,0)[rb]{\DeformationFile[scale=0.7]{02}{0001}{17}}}
\put(-0.01,3.05){\makebox(0,0)[rb]{\DeformationFile[scale=0.7]{02}{0001}{22}}}
\put(-0.01,2.05){\makebox(0,0)[rb]{\DeformationFile[scale=0.7]{02}{0001}{28}}}
\put(-0.01,1.05){\makebox(0,0)[rb]{\DeformationFile[scale=0.7]{02}{0001}{32}}}
\put(0.01,5.05){\makebox(0,0)[lb]{\DeformationFile[scale=0.7]{11}{0001}{3}}}
\put(0.01,4.05){\makebox(0,0)[lb]{\DeformationFile[scale=0.7]{11}{0001}{5}}}
\put(0.01,3.05){\makebox(0,0)[lb]{\DeformationFile[scale=0.7]{11}{0001}{11}}}
\put(0.01,2.05){\makebox(0,0)[lb]{\DeformationFile[scale=0.7]{11}{0001}{36}}}
\put(0.01,1.05){\makebox(0,0)[lb]{\DeformationFile[scale=0.7]{11}{0001}{47}}}
\put(0.01,0.05){\makebox(0,0)[lb]{\DeformationFile[scale=0.7]{11}{0001}{109}}}
\put(-.8,-.1){\makebox(0,0)[tr]{(a) small friction $\mu{=}0.2$}}
\put(.8,-.1){\makebox(0,0)[tl]{(b) large friction $\mu{=}1.1$}}
\end{picture}
\vspace*{.5em}
\caption{Deformed configurations of the skewed blocks and contact traction distributions at various instants of the loading process,   smaller viscosity $\chi{=}0.0001\,$s.
We can see not only stick/slide regimes but also jump on the last snapshots for large friction.
}\label{Fig_SkewDefMU}
\end{figure}
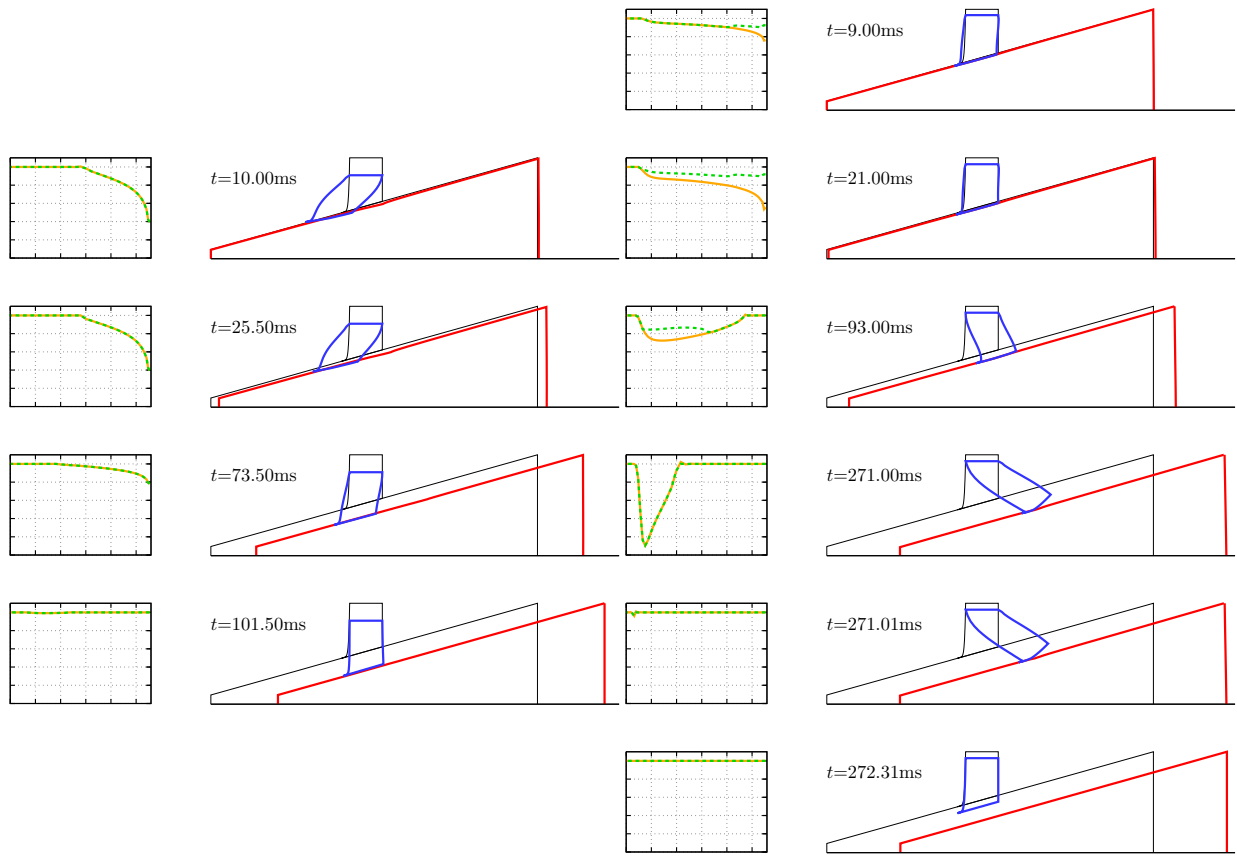

Figure~\ref{Fig_SkewDefMU} shows  the deformed configurations of the punch and a top part of the bottom block.
It should be noted that, unlike the previous examples, where the time step was fixed, here the adaptive procedure in time is used.
The graphs show  various stages in the evolution of  deformation and qualitatively also of the contact traction distribution for $\mu{=}0.2$ and for $\mu{=}1.1$.

The displacements for the small-friction case in Figure~\ref{Fig_SkewDefMU}(a) are scaled by a factor of 1000 to see what happens.
This factor also causes that the bodies seem to be slightly interpenetrated in some time steps.
The graphs of tractions contain both tangential and $\mu$-scaled normal components.
As there is only one curve visible during the whole load history, the bodies slides on each other.

A magnifying factor for the displacements in the large-friction case in Figure~\ref{Fig_SkewDefMU}(b) is 350.
The graphs of tractions now document that in the beginning of the loading there is a sliding contact zone close to the nose $N$ (this was the reason for including it) and a stress singularity at the other end point of the contact zone.
When the bottom block is pushed, the stresses at this end point decreases until the right part of the punch detaches.
There remains only a stress concentration near to the nose end.
Nevertheless, the two bodies remain in mutual contact until the tangential component of the traction reaches its sliding threshold and then they separate.
Let us pay attention especially to the last three used time instants.
The very last one corresponds to the instant of returning back to the original shape, due to (small) viscosity.
The other two instants, in fact very close to each other, pertain to the abrupt jump of the total force as seen e.g. in Figure~\ref{Fig_SkewEFTR}.
The distribution of the contact traction also abruptly changes to zero.
We stress once more that this is a special behaviour of the contact associated only to large friction.

\begin{figure}[ht]
\centering
\PictureFile{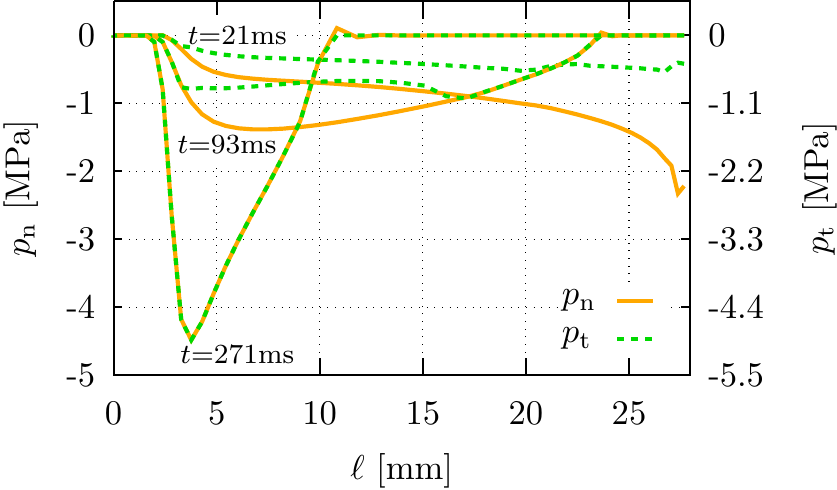}
\caption{A detailed comparison of contact tractions at three selected instants from Fig.~\ref{Fig_SkewDefMU}(b).
}\label{Fig_SkewContact}
\end{figure}

After having seen qualitative distributions of the contact tractions during the loading histories, it is worth to present also their quantitative representations. 
The distributions of contact tractions  for the more interesting case $\mu{=}1.1$ with the largest variations at three instants are shown in Figure~\ref{Fig_SkewContact}.
In the first plotted time step $t{=}21$ms which is close to the initiation of horizontal pushing, there is only a small area of the slip contact.
Chronologically, the second time step $t{=}93$ms presents extension of the slip zone and also detachment of the right part of the interface, where the contact tractions vanish.
The last shown time step $t{=}271$ms corresponds to a situation before the total detachment of the blocks and the overall state of the contact is slipping.

\section{Conclusions}

An energy based model for solving frictional contact problem has been proposed, analysed and implemented in a computational code.
The model uses a regularization of  classical contact conditions by allowing a small interpenetration applying a normal-compliance model and also considers a visco-elastic material to make the solution more regular.
The numerical implementation of spatial discretization {\em via} SGBEM has permitted the whole problem to be defined and solved  only by working with the boundary and contact zone data.
Additionally, some algebraic manipulation have allowed us to use SGBEM for an elasto-static problem and QP.
Two examples of various contact configurations, some of them quite difficult and intricate, have been used to validate  the model and its computational implementation.
It is remarkable that the proposed model provides satisfactory results also for large friction coefficients, even greater than one, which was documented in the computational analysis as well and that the problems with such large friction coefficients can lead to totally different results than those with small ones.
Finally, one could also take mechanical damaging due to micro-cracks at the interface into account, similarly as in e.g.~\cite{HaShSo01VNAQ}, which will be discussed in a more general fashion in the forthcoming paper~\cite{VoMaRo??}.

\subsubsection*{Acknowledgement}
This research has been covered by the grant Ministry of Education of the Slovak Republic (VEGA 1/0477/15 and 1/0078/16), by the University of Seville (SAB2010-0082),
by the Junta de Andaluc\'{\i}a and European Social Fund (TEP-4051), the Spanish Ministry of Economy and Competitiveness and   European Regional Development Fund  (MAT2012-37387, MAT2015-71036-P), and by the Czech Science Foundation through the grant 16-03823.
Both R.~V.\ and T.~R.\ acknowledge hospitality of the University of Seville during their stays in 2011-2015.
Beside, T.~R.\ acknowledges also the institutional support RVO: 61388998 (\v CR).

\bibliography{VodickaManticRoubicek_CAM_final}

\appendix

\section{Implementation of SGBEM}\label{App_SGBEM}

For the solution with respect to bulk domains, the SGBEM~\cite{bonnet98A1,sutradhar08B1} implementation, deduced from the energetic principles as shown in~\cite{vodicka07A1,vodicka11A1} and guaranteeing the positive-definite character of the computed strain energy is used.
Let us briefly summarized some details required for such implementation.

The BIEs solved by SGBEM are the {\em Somigliana displacement and traction identities}, written for each particular domain $\Omega^\eta$ separately:
\begin{subequations}\label{Eq_utequations}
\begin{align}
\label{Eq_ueq}
 \frac12 v^\eta_k(x) & =\int_{\Gamma^\eta}\!\!U^\eta_{kl}(x,y) \P^\eta_l(y)\dd\Gamma(y)
- \cint{\Gamma^\eta}\!\!T^\eta_{kl}(x,y)v^\eta_l(y)\dd\Gamma(y),&
\text{for a.a. }x&\!\in\!
\GD^\eta\cup\GC,
\\\label{Eq_teq}
\frac12 \P^\eta_k(x)&=\cint{\Gamma^\eta}\!\!T^{\eta\ast}_{kl}(x,y)\P^\eta_l(y)\dd\Gamma(y)
- \fint{\Gamma^\eta}\!\!S^\eta_{kl}(x,y)v^\eta_l(y)\dd\Gamma(y),&
\text{for a.a. }x&\!\in\!\GN^\eta.
 \end{align}
\end{subequations}
The integral kernels in the above equations are  the (weakly singular) Kelvin fundamental solution $U^\eta_{ij}(x,y)$ -- the response of the elastic plane to a point load, and the associated derivative kernels obtained by the differential traction operator applied with respect to one or both variables -- the strongly singular kernel $T^\eta_{ij}(x,y)$ and the hypersingular kernel $S^\eta_{ij}(x,y)$.
Due to the integral kernel singularities, the integral denoted by $\cint{\Gamma}$ or ${\fint{\Gamma^{}}}$, respectively, stands for the Cauchy principal value or the Hadamard finite part of the integral.
It should also be noted that the formulation of BIE may provide non-unique solutions for a uniquely solvable BVP if a domain contains a cavity with prescribed Neumann boundary conditions.
Such non-uniqueness can be removed by techniques derived in~\cite{vodicka06A1}.

Introducing the operator notation will allow us to rewrite the weighted formulation of the BIE system in a compact and transparent form.
Let
\begin{equation}\label{Eq_FormalOperator}
\VectorSymbol{\omega}^{\eta \top}_q\Vector{Z}^\eta_{qr}\Vector{\DIR}^\eta_r=
\int_{\Gamma^\eta_q}\omega^\eta_j(y)\left(\int_{\Gamma^\eta_r}Z^\eta_{ji}(y,x)w^\eta_i(x)\dd\Gamma(x)
\right)\dd\Gamma(y),
\end{equation}
where $\omega$ stands for $\varphi$ or $\psi$, while $w$ stands for $v$ or $p$, and further $q$ and $r$ stand for D, N, and C, and eventually $Z^\eta$ stands for $U^\eta$, $T^\eta$, $T^{\eta\ast}$ or $S^\eta$, and where the inner integral can be regular, weakly singular, Cauchy principal value or Hadamard finite part integral.

The present system of BIEs written in the weighted formulation can be arranged in the following block form, cf.~\cite{vodicka07A1,vodicka11A1}:
\begin{equation}\label{RepresentationMatrix}
\begin{pmatrix}
\VectorSymbol{\varphi}\AA\D \\ \VectorSymbol{\psi}\AA\N \\ \VectorSymbol{\varphi}\AA\C \\ \VectorSymbol{\psi}\AA\C\\
\VectorSymbol{\varphi}\BB\D \\ \VectorSymbol{\psi}\BB\N \\ \VectorSymbol{\varphi}\BB\C \\ \VectorSymbol{\psi}\BB\C
\end{pmatrix}^{\negthickspace\top}\negthickspace\negthickspace
\begin{pmatrix}
\Vector{K}\AA & \Vector{M}\AB \\
\Vector{M}\ABast & \Vector{K}\BB
\end{pmatrix}
\begin{pmatrix}
\Vector{\P}\AA\D \\ \Vector{v}\AA\N \\ \Vector{\P}\AA\C \\ \Vector{v}\AA\C\\
\Vector{\P}\BB\D \\ \Vector{v}\BB\N \\ \Vector{\P}\BB\C \\ \Vector{v}\BB\C
\end{pmatrix}
=
\begin{pmatrix}
\VectorSymbol{\varphi}\AA\D \\ \VectorSymbol{\psi}\AA\N \\ \VectorSymbol{\varphi}\AA\C \\ \VectorSymbol{\psi}\AA\C\\
\VectorSymbol{\varphi}\BB\D \\ \VectorSymbol{\psi}\BB\N \\ \VectorSymbol{\varphi}\BB\C \\ \VectorSymbol{\psi}\BB\C
\end{pmatrix}^\top
\begin{pmatrix}
\Vector{H}\AA & \Vector{0} & -\left(\Vector{M}\AB\right)_{\cdot4} \\
\Vector{0} & \Vector{H}\BB & \Vector{0}
\end{pmatrix}
\begin{pmatrix}
\wt{\Vector{\DIR}}_\Dir\AA \\ \Vector{f}\N\AA\\ \wt{\Vector{\DIR}}_\Dir\BB \\ \Vector{f}\N\BB\\ \Vector{\W}
\end{pmatrix},
\end{equation}
with
\begin{gather}\label{Eq_ReprMatricesDef}
\Vector{K}^\eta{=}
\begin{pmatrix}
-\Vector{U}^\eta\DD & \Vector{T}^\eta\DN & -\Vector{U}^\eta\DC & \Vector{T}^\eta\DC\\
 \Vector{T}^{\eta \ast}\ND & -\Vector{S}^\eta\NN & \Vector{T}^{\eta \ast}\NC & -\Vector{S}^\eta\NC\\
-\Vector{U}^\eta\CD &\Vector{T}^\eta\CN & -\Vector{U}^\eta\CC &\omega^\eta\frac12\Vector{I}^\eta\CC{+}\Vector{T}^\eta\CC\\
 \Vector{T}^{\eta \ast}\CD & -\Vector{S}^\eta\CN & \omega^\eta\frac12\Vector{I}^\eta\CC{+}\Vector{T}^{\eta \ast}\CC & -\Vector{S}^\eta\CC
\end{pmatrix},\ \ \
\omega^\eta=\begin{cases}\ -1&\text{if }\eta=A,\\[-.3em] 1&\text{if }\eta=B,\end{cases}
 \nonumber\\[-1em]
{} \\
\Vector{H}^\eta{=}
\begin{pmatrix}
 -\frac12\Vector{I}^\eta\DD{-}\Vector{T}^\eta\DD & \Vector{U}^\eta\DN \\
 \Vector{S}^\eta\DD&  \frac12\Vector{I}^\eta\NN{-}\Vector{T}^{\eta \ast}\NN \\
 -\Vector{T}^\eta\CD& \Vector{U}^\eta\CN \\
 \Vector{S}^\eta\CD&  -\Vector{T}^{\eta \ast}\CN
\end{pmatrix}, \qquad\qquad\quad
\Vector{M}\AB{=}
\begin{pmatrix}
\Vector{0} & \Vector{0} & \Vector{0} & \Vector{0}\\
\Vector{0} & \Vector{0} & \Vector{0} & \Vector{0}\\
\Vector{0} & \Vector{0} & \Vector{0} & \Vector{I}\AB\CC\\
\Vector{0} & \Vector{0} & \Vector{0} & \Vector{0}
\end{pmatrix}.\nonumber\qquad
\end{gather}
In the previous relations $\Vector{I}^\eta$ denotes the identity operator with the subscripts and superscripts specifying the part of the boundary where it is restricted.
The novelty of the BIE system in \eqref{RepresentationMatrix} with respect to that developed in~\cite{vodicka07A1,vodicka11A1}  is that no displacement gap at $\GC$ was considered therein.
The present formulation, unlike that of~\cite{PaMaRo14BAMI}, optimizes the energy in terms of the displacement gap $\Vector{\W}$, so that the number of unknowns in the optimization process is roughly a half of those used in the cited reference.
On the other hand, in the solution of the elastic BVP a multi-domain approach is used in the present formulation, while the aforementioned reference solves the BVP for each domain separately.
The BIE system~\eqref{RepresentationMatrix} will be solved numerically by SGBEM.
To this end,    the variables appearing there are approximated by linear continuous boundary elements~\cite{paris97} (allowing discontinuities of the tractions at the
junctions of the elements if required).
The approximation formulas can be written in the form
\begin{equation}\label{Eq_UTapproximation}
\Vector{v}^\eta(x)=\sum_n\Vector{N}^\eta_{\psi n}(x)\VectorMatrix{v}^\eta_n,\quad
\Vector{\W}(x)=\sum_{n\C^{}}\Vector{N}\BB_{\psi {n\C^{}}}(x)\VectorMatrix{\W}_{n\C^{}},\quad
\Vector{\P}^\eta(x)=\sum_m\Vector{N}^\eta_{\varphi m}(x)\VectorMatrix{\P}^\eta_m,
\end{equation}
where $\Vector{N}^\eta_{\psi\ell}(x)$ and $\Vector{N}^\eta_{\varphi\ell}(x)$, respectively, are matrices containing the shape functions of displacements and tractions associated to node $\ell$ at $x_\ell^\eta{\in}\Gamma^\eta$, and $\VectorMatrix{v}^\eta_\ell$, $\VectorMatrix{\W}_{\ell}$ and $\VectorMatrix{\P}^\eta_\ell$, respectively, are vectors containing the components of the displacement, displacement gap and traction vectors at the node $\ell$.
Let $\VectorMatrix{v}^\eta$, $\VectorMatrix{\DIR}_\Dir^\eta$, $\VectorMatrix{\P}^\eta$, $\VectorMatrix{f}\N^\eta$ and $\VectorMatrix{h}$, respectively, denote the vectors
containing all unknown nodal displacements (transformed by~\eqref{Eq_SubstV}), all prescribed nodal  displacements (transformed as in~\eqref{Eq_AdmissibleV}), all unknown nodal tractions, all prescribed nodal tractions associated to $\Gamma^\eta$ and all fictitious nodal displacements gaps at $\GC$.
Let the subvectors of the nodal unknowns at the boundary parts $\GD^\eta$, $\GN^\eta$, and $\GC$, respectively, be distinguished by the same subscripts D, N, and C, respectively.
The set of vectors of virtual functions $\VectorSymbol{\psi}^\eta$ and $\VectorSymbol{\varphi}^\eta$ can be chosen in the way that they are equal to the shape functions associated to each nodal unknown.
Such a choice leads to the following system of linear algebraic equations with a  symmetric  matrix:
\begin{equation}\label{DiscretizationMatrix}
\begin{pmatrix}
\Matrix{K}\AA & \Matrix{M}\AB \\
\Matrix{M}^{AB \ast} & \Matrix{K}\BB
\end{pmatrix}
\begin{pmatrix}
\VectorMatrix{\P}\AA\D \\[.1em] \VectorMatrix{v}\AA\N \\[.1em] \VectorMatrix{\P}\AA\C \\[.1em] \VectorMatrix{v}\AA\C\\[.1em]
\VectorMatrix{\P}\BB\D \\[.1em] \VectorMatrix{v}\BB\N \\[.1em] \VectorMatrix{\P}\BB\C \\[.1em] \VectorMatrix{v}\BB\C
\end{pmatrix}
=
\begin{pmatrix}
\Matrix{H}\AA & \Matrix{0} & -\left(\Matrix{M}\AB\right)_{\cdot4} \\
\Matrix{0} & \Matrix{H}\BB & \Matrix{0}
\end{pmatrix}
\begin{pmatrix}
\wt{\VectorMatrix{\DIR}}_\Dir\AA \\ \VectorMatrix{f}\N\AA\\ \wt{\VectorMatrix{\DIR}}_\Dir\BB \\ \VectorMatrix{f}\N\BB\\ \VectorMatrix{\W}
\end{pmatrix},
\end{equation}
where
\begin{gather}\label{Eq_ReprMatricesDefDiscr}
\Matrix{K}^\eta{=}
\begin{pmatrix}
-\Matrix{U}^\eta\DD & \Matrix{T}^\eta\DN & -\Matrix{U}^\eta\DC\ \  & \Matrix{T}^\eta\DC\\
 \Matrix{T}^{\eta \ast}\ND & -\Matrix{S}^\eta\NN\ \  & \Matrix{T}^{\eta \ast}\NC & -\Matrix{S}^\eta\NC\ \ \\
-\Matrix{U}^\eta\CD &\Matrix{T}^\eta\CN & -\Matrix{U}^\eta\CC\ \  &\omega^\eta\frac12\Matrix{M}^\eta\CC{+}\Matrix{T}^\eta\CC\\
 \Matrix{T}^{\eta \ast}\CD & -\Matrix{S}^\eta\CN\ \  & \omega^\eta\frac12\Matrix{M}^\eta\CC{+}\Matrix{T}^{\eta \ast}\CC & -\Matrix{S}^\eta\CC\ \
\end{pmatrix},\ \ \ \ \ \eta={\rm A,B}, \nonumber\\[-1em]
{} \\
\Matrix{H}^\eta{=}
\begin{pmatrix}
 -\frac12\Matrix{M}^\eta\DD{-}\Matrix{T}^\eta\DD & \Matrix{U}^\eta\DN \\
 \Matrix{S}^\eta\DD&  \frac12\Matrix{M}^\eta\NN{-}\Matrix{T}^{\eta \ast}\NN \\
 -\Matrix{T}^\eta\CD\ \ & \Matrix{U}^\eta\CN \\
 \Matrix{S}^\eta\CD&  -\Matrix{T}^{\eta \ast}\DN\ \
\end{pmatrix}, \qquad
\Matrix{M}\AB{=}
\begin{pmatrix}
\Matrix{0} & \Matrix{0} & \Matrix{0} & \Matrix{0}\\
\Matrix{0} & \Matrix{0} & \Matrix{0} & \Matrix{0}\\
\Matrix{0} & \Matrix{0} & \Matrix{0} & \Matrix{M}\AB\CC\\
\Matrix{0} & \Matrix{0} & \Matrix{0} & \Matrix{0}
\end{pmatrix}.\nonumber
\end{gather}

The elements of the submatrices denoted with letters \Matrix{U}, \Matrix{T} and \Matrix{S} are formed by double integrals including the integral kernel denoted by the same letter as is usual in SGBEM, see also~\eqref{Eq_FormalOperator}.
The square $d{\times}d$  submatrices, associated to nodes $n$ and $m$, of the  mass matrices $\Matrix{M}_{rr}$, with $r$ being D, N, or C, are defined by the integrals:
\begin{equation}\label{MatrixMeta}
(\Matrix{M}^\eta_{rr})_{mn}^{}=
\int_{\Gamma^\eta_r}\Vector{N}^\eta_{\varphi m}(x)\Vector{N}^\eta_{\psi n}(x)
\,\dd\Gamma, \quad
\qquad
(\Matrix{M}\AB\CC)_{mn}^{}=
\int_{\GC}\Vector{N}\AA_{\varphi m}(x)\Vector{N}\BB_{\psi n}(x)\,\dd\Gamma.
\end{equation}

\def\V{{\tilde{v}}}

\section{An expression for a generalized Poincar\'e-Steklov operator $\calPt$}\label{Sec_PCC}

First, let us consider BVPs with boundary conditions~\eqref{BVP-Dirichlet} and~\eqref{BVP-Neuman} in both domains $\Omega^\eta$ independently considering $\GC$ as a Dirichlet boundary part in addition to $\GD^\eta$.
The pertinent Poincar\'e-Steklov operators considered as invertible operators $\calP^\eta$ from a quotient space of $H^{1/2}( \GC\cup \GD^\eta\cup \GN^\eta; \R^d)$ by the space of rigid body motions to the space $H^{-1/2}( \GC\cup \GD^\eta\cup \GN^\eta; \R^d)$ can be represented in  a block form, cf.~\cite{Agranovich2011,Helsing2011},   see also \cite{PaMaRo13CM},
\begin{equation}\label{Eq_PSblock}
\begin{pmatrix}\P^\eta_{\D}\\[.2em] \P^\eta{\N} \\[.2em] \P^\eta{\C}\end{pmatrix}=
\begin{pmatrix}
\calP^\eta_{\DD}\!&\!\calP^\eta_{\DN}\!&\!\calP^\eta_{\DC}
\\[.2em] \calP^\eta_{\ND}\!&\!\calP^\eta_{\NN}\!&\!\calP^\eta_{\NC}
\\[.2em] \calP^\eta_{\CD}\!&\!\calP^\eta_{\CN}\!&\!\calP^\eta_{\CC}
\end{pmatrix}
\begin{pmatrix}v^\eta_{\D}\\[.2em] v^\eta_{\N}\\[.2em] v^\eta_{\C}\end{pmatrix}.
\end{equation}
The above operator can be used to map prescribed boundary data ${\Vector{\DIR}}^\eta_\Dir$ at $\GD^\eta$, $f\N^\eta$ at $\GN^\eta$ and $v_{\C}$ at $\GC$ to the unknown ones.
The operator meant as a generalization of the Poincar\'e-Steklov operator is denoted $\calPt^\eta$ and its block structure can formally be obtained from~\eqref{Eq_PSblock} as
\begin{multline}\label{Eq_Sblock}
\begin{pmatrix}\P^\eta{\D}\\ v^\eta{\N} \\[.2em] \P^\eta{\C}\end{pmatrix}=
\begin{pmatrix}
\calPt^\eta_{\DD}\!&\!\calPt^\eta_{\DN}\!&\!\calPt^\eta_{\DC}\\[.2em] \calPt^\eta_{\ND}\!&\!\calPt^\eta_{\NN}\!&\!\calPt^\eta_{\NC}\\[.2em] \calPt^\eta_{\CD}\!&\!\calPt^\eta_{\CN}\!&\!\calPt^\eta_{\CC}
\end{pmatrix}
\begin{pmatrix}{\Vector{\DIR}}^\eta_\Dir\\[.2em] f\N^\eta \\[.2em] v^\eta_{\C}\end{pmatrix}\\
=
\begin{pmatrix}
\calP^\eta_{\DD}{-}{\calP^\eta_{\DN}}{\left(\calP^\eta_{\NN}\right)^{-1}}\!{\calP^\eta_{\ND}} & \calP^\eta_{\DN}{\left(\calP^\eta_{\NN}\right)^{-1}} & \calP^\eta_{\DC}{-}{\calP^\eta_{\DN}}{\left(\calP^\eta_{\NN}\right)^{-1}}\!{\calP^\eta_{\NC}}\\[.2em]
{\left(\calP^\eta_{\NN}\right)^{-1}}\!{\calP^\eta_{\ND}} & {\left(\calP^\eta_{\NN}\right)^{-1}} & {\left(\calP^\eta_{\NN}\right)^{-1}}\!\calP^\eta_{\NC}\\[.2em]
\calP^\eta_{\CD}{-}{\calP^\eta_{\CN}}{\left(\calP^\eta_{\NN}\right)^{-1}}\!{\calP^\eta_{\ND}} & \calP^\eta_{\CN}{\left(\calP^\eta_{\NN}\right)^{-1}} & \calP^\eta_{\CC}{-}{\calP^\eta_{\CN}}{\left(\calP^\eta_{\NN}\right)^{-1}}\!{\calP^\eta_{\NC}}
\end{pmatrix}
\begin{pmatrix}{\Vector{\DIR}}^\eta_\Dir\\[.2em] f\N^\eta \\[.2em] v^\eta_{\C}\end{pmatrix}.
\end{multline}
Knowing the mapping properties and symmetry of the Poincar\'e-Steklov restrictions, which are, e.g.,
$\calP^\eta\CD={\calP^\eta\DC}^*:H^{1/2}(\GD^\eta;\R^d) \mapsto H^{-1/2}(\GC;\R^d)$, we can recognize the symmetry of $\calPt^\eta$ considered as the following map:
\begin{multline}\label{Eq_Smap}
\calPt^\eta:H^{1/2}(\GD^\eta;\R^d){\times}H^{-1/2}(\GN^\eta;\R^d){\times}H^{1/2}(\GC;\R^d)\\ \mapsto H^{-1/2}(\GD^\eta;\R^d){\times}H^{1/2}(\GN^\eta;\R^d){\times}H^{-1/2}(\GC;\R^d).
\end{multline}
The restriction to $\GC$ is denoted as $\calPt_{\CC}^\eta$. It should be noted
that it requires no restriction to a quotient space of  $H^{1/2}( \GC;\R^d)$
because the functions are considered vanishing in nonempty $\GD^\eta$.

Now, let us consider the TBVP~\eqref{BVP-k}.
The generalized Poincar\'e-Steklov operator $\calPt$ of the same block structure as defined in  \eqref{Eq_Sblock} can be introduced.
Having in mind the equilibrium condition $\P\AA\C{+}\P\BB\C{=}0$ and the interface gap $\Vector{\W}$: $v\AA\C{-}v\BB\C{=}\Vector{\W}$  of~\eqref{BVP-k} and introducing $\V{=}v\AA\C{+}v\BB\C$ and $\P\C{=}\P\AA\C{=}{-}\P\BB\C$, we obtain
\begin{align}\label{Eq_PS-blockTwo}
\begin{pmatrix}\left(\begin{smallmatrix}  \P\D\AA \\  v\AA\N \\  \P\C\end{smallmatrix}\right)
\\[.9em]
\left(\begin{smallmatrix} \P\BB\D \\  v\BB\N \\-\P\C\end{smallmatrix}\right)
\end{pmatrix}
=
\begin{pmatrix}
\calPt\AA \!&\! 0   \\
0 \!&\! \calPt\BB
\end{pmatrix}
\begin{pmatrix}
\left(\begin{smallmatrix}{\Vector{\DIR}}\AA_\Dir \\ f\N\AA \\ \frac{\W+\V}2\end{smallmatrix}\right) \\[.9em]
\left(\begin{smallmatrix}{\Vector{\DIR}}\BB_\Dir \\  f\N\BB \\ \frac{-\W+\V}2 \end{smallmatrix}\right)
\end{pmatrix}.
\end{align}
Then, the block equation can be rearranged  to
\begin{align}\label{Eq_PS-blockTwovar}\setlength{\arraycolsep}{1pt}
\begin{pmatrix}  \P\D\AA \\[.5em]  \P\BB\D \\[.5em] v\AA\N \\[.5em] v\BB\N \\[.5em]  \P\C \\[.5em]   0\end{pmatrix}{=}
\begin{pmatrix}
\calPt\AA\DD & 0 & \calPt\AA\DN & 0 &\frac12\calPt\AA\DC& \frac12\calPt\AA\DC\\
0 & \calPt\BB_{\DD} & 0 & \calPt\BB_{\DN} & -\frac12\calPt\BB_{\DC}& \frac12\calPt\BB_{\DC}\\
\calPt\AA\ND & 0 & \calPt\AA\NN & 0 &\frac12\calPt\AA\NC& \frac12\calPt\AA\NC\\
0 & \calPt\BB_{\ND} & 0 & \calPt\BB_{\NN} & -\frac12\calPt\BB_{\NC}& \frac12\calPt\BB_{\NC}\\
\frac12\calPt\AA_{\CD} & -\frac12\calPt\BB_{\CD} & \frac12\calPt\AA_{\CN} & -\frac12\calPt\BB_{\CN} & \frac14\!\left(\calPt\AA_{\CC}{+}\calPt\BB_{\CC}\right) & \frac14\!\left(\calPt\AA_{\CC}{-}\calPt\BB_{\CC}\right) \\
\frac12\calPt\AA_{\CD} & \frac12\calPt\BB_{\CD} & \frac12\calPt\AA_{\CN} & \frac12\calPt\BB_{\CN} & \frac14\!\left(\calPt\AA_{\CC}{-}\calPt\BB_{\CC}\right) & \frac14\!\left(\calPt\AA_{\CC}{+}\calPt\BB_{\CC}\right)
\end{pmatrix}\!
\begin{pmatrix}
{\Vector{\DIR}}\AA_\Dir \\[.4em] {\Vector{\DIR}}\BB_\Dir \\[.4em] f\N\AA \\[.4em] f\N\BB \\[.4em] \frac{\W+\V}2 \\[.4em] \frac{-\W+\V}2
\end{pmatrix}.
\end{align}
The last row and column can be eliminated due to invertibility of $\frac14\left(\calP\AA\CC{+}\calP\BB\CC\right)$ so that the resulting expression is similar to \eqref{Eq_Sblock},
\begin{equation}\label{Eq_SblockTwo}
\begin{pmatrix}\P{\D}\\[.3em] v{\N} \\[.3em] \P{\C}\end{pmatrix}{=}
\begin{pmatrix}
\calPt_{\DD}\!&\!\calPt_{\DN}\!&\!\calPt_{\DC}\\ \calPt_{\ND}\!&\!\calPt_{\NN}\!&\!\calPt_{\NC}\\ \calPt_{\CD}\!&\!\calPt_{\CN}\!&\!\calPt_{\CC}
\end{pmatrix}\negthickspace\begin{pmatrix}{\Vector{\DIR}}_\Dir\\[.3em] f\N \\[.3em] \W\C\end{pmatrix}
=\calPt
\begin{pmatrix}{\Vector{\DIR}}_\Dir\\[.3em] f\N \\[.3em] \W\C\end{pmatrix}
\quad\text{ with }  \begin{alignedat}{2} \P{\D}&=\begin{pmatrix} \P\AA\D\\ \P\BB\D\end{pmatrix}, &\ v\N&=\begin{pmatrix} v\AA\N \\ v\BB\N\end{pmatrix}, \\
f\N&=\begin{pmatrix} f\AA\N \\ f\BB\N\end{pmatrix}, & \ {\Vector{\DIR}}_\Dir&=\begin{pmatrix} {\Vector{\DIR}}\AA_\Dir \\ {\Vector{\DIR}}\BB_\Dir\end{pmatrix}.
\end{alignedat}
\end{equation}
Note that the eliminated last row of the left-hand side was zero so it has no influence on the rest of the left-hand side in such an elimination and that $\V$ does not appear in the resulting expression.
In fact, this operator $\calPt$ is used in \eqref{Eq_ChangeOmegaToGamma} and referred to in the lines below.
The above elimination produces $\calPt_{\CC}$ in the form:
\begin{equation}\label{Eq_Sdef}
\calPt\CC=\frac14\left[ \left(\calPt\CC\AA + \calPt\CC\BB\right) - \left(\calPt\CC\AA - \calPt\CC\BB\right)\left(\calPt\CC\AA + \calPt\CC\BB\right)^{-1}\left(\calPt\CC\AA -\calPt\CC\BB\right) \right].
\end{equation}

This restriction of the generalized Poincar\'e-Steklov operator $\calPt$  to $H^{1/2}(\GC;\R^d)$ can also be introduced in another way, assuming a simplified case of $g\D^\eta{=}0$ and $f\N^\eta{=}0$.
Having in mind the  conditions $v\AA\C{-}v\BB\C{=}\Vector{\W}$,  $\P\AA\C{+}\P\BB\C{=}0$ and the restricted generalized Poincar\'e-Steklov operators to $\GC$ which are defined due to  \eqref{Eq_Sblock} as $\calPt_{\CC}\AA(v\AA\C){=}\P\AA\C$ and $\calPt_{\CC}\AA(v\BB\C){=}\P\BB\C$, we obtain
\begin{multline}\label{Eq_PSmulti}
\Vector{\W}= v\AA\C-v\BB\C= \left(\calPt_{\CC}\AA\right)^{-1}\left(\P\AA\C\right) - \left(\calPt_{\CC}\BB\right)^{-1}\left(\P\BB\C\right)=
\left(\calPt_{\CC}\AA\right)^{-1}\left(\P\AA\C\right) + \left(\calPt_{\CC}\BB\right)^{-1}\left(\P\AA\C\right)=\\
\left[ \left(\calPt_{\CC}\AA\right)^{-1} + \left(\calPt_{\CC}\BB\right)^{-1} \right]\left(\P\AA\C\right)
=-\left[ \left(\calPt_{\CC}\AA\right)^{-1} + \left(\calPt_{\CC}\BB\right)^{-1} \right]\left(\P\BB\C\right).
\end{multline}
Therefore, we can write
$[(\calPt_{\CC}\AA)^{-1}{+}(\calPt_{\CC}\BB)^{-1}]^{-1}\Vector{\W}=\P=\P\AA\C={-}\P\BB\C$.
This relation leads to another equivalent representation of the operator defined in \eqref{Eq_Sdef},
\begin{equation}\label{Eq_PSdef}
\calPt_{\CC}=\left[ \left(\calPt_{\CC}\AA\right)^{-1}\!+ \left(\calPt_{\CC}\BB\right)^{-1} \right]^{-1}.
\end{equation}
Here, the operator is defined as a double of the harmonic mean of the restricted operators associated to each particular subdomain $\calPt_{\CC}^\eta$.
The properties of these operators, e.g.  coercivity or symmetry, which are required for a correct treatment of the functional in \eqref{Eq_ChangeOmegaToGamma}, are then inherited also by the operator  $\calPt_{\CC}$.

Note also that the solution at $\GC$ corresponding  to each particular subdomain, can be expressed as:
\begin{equation}\label{Eq_vSmap}
v\AA\C= \left(\calPt_{\CC}\AA\right)^{-1}\!\calPt_{\CC}\left(\Vector{\W}\right),
\qquad
v\BB\C=-\left(\calPt_{\CC}\BB\right)^{-1}\!\calPt_{\CC}\left(\Vector{\W}\right)
\end{equation}
in the above simplified case.

Similar ideas can be used also for the discretized versions of the operator $\calPt$ in relation to the discretized operator $\calPt_h$ introduced in~\eqref{def-of-K-h} and the lines below.
As long as in each particular subdomain we obtain, using appropriate discretizations, uniformly positive approximations of the operators say $\calPt_{\CC}^\eta$, also the discretization operators $\left({\calPt_h}\right){\CC}$ can be constructed uniformly positive.

\section{Theoretical justification: notes on numerical stability and convergence}\label{Sec_Theory}

Let us briefly overview some theoretical aspects of the model and its numerical implementation.
Now, numerical solutions with different  time step $\tau>0$ as well as space-discretization parameter $h>0$ will be considered.
We denote by $\baru_{\tau h}$,  $\barz_{\tau h}$, and $\barw_{\tau h}$ the piecewise constant interpolants which use the values of the numerical solution at superior times of the pertinent time step intervals, and we also define  $\underline \Z_{\tau h}$ as the piecewise constant interpolant which uses the value $\Z^{k-1}_{\tau h}$ at inferior times of the time step interval $((k{-}1)\tau,k\tau)$. Occasionally, we also use the notation $u_{\tau h}$
or $\Z_{\tau h}$ for the continuous, piecewise affine interpolants.
Similarly, we will define $\barf_\tau$, $\barF_\TAU^{}$, etc..
We will consider a general dimension $d{=}2,3$ in this section.

First, from the energy estimate \eqref{Eq_EbilanceDiscr}, which holds for the approximate solution too, one gets the a-priori estimate for the sequence $\{u_{\tau h}\}_{\tau>0,h>0}$, which can be proved bounded in the space $H^1(0,T;H^1(\Omega\AA\cup\Omega\BB;\R^d))$ provided that ${\DIR}_\Dir\in H^1(0,T;H^{1/2}(\GD;\R^d))$ and $f\N\in H^1(0,T;H^{-1/2}(\GN;\R^d))$.

The weak convergence of the time-discretization itself is quite simple because
the nonlinear terms on $\GC$ are of lower order thanks to the regularizing
character of the normal-compliance model. In fact, a strong convergence in
$H^1(0,T;H^1(\Omega\AA\cup\Omega\BB;\R^d))$ can be derived, cf.~also
\cite{HaShSo01VNAQ} where, under zero-Dirichlet condition, a similar
semi-implicit time discretization (actually combined with a FEM-discretization)
has been used and rate of convergence towards the unique solution has been obtained. Here, an important message from \cite{HaShSo01VNAQ} is that, in particular, the original problem \eqref{BVP} has a unique weak solution under the mentioned natural data qualification; in fact, as \cite{HaShSo01VNAQ} uses ${\DIR}_\Dir{=}0$, the uniqueness of the weak solution of the problem \eqref{BVP} can be seen from \cite{HaShSo01VNAQ} only after transformation to the zero-Dirichlet condition.

The convergence proof of the above SGBEM-approximation is more complicated than the FEM-approximation in \cite{HaShSo01VNAQ} because \eqref{Eq_ChangeOmegaToGamma} and thus also \eqref{Eq_DiscreteGeneral} are only approximations of \eqref{Eq_DiscreteFunctional}, in contrast to the conformal FEM without numerical integration as used in \cite{HaShSo01VNAQ}.
Also, in order to use the static BEM, we made several transformations in Section~\ref{SecImplTime-BIE} and must now return back to the form which the original problem has. To this goal, let us first write \eqref{Biot-disc-g-h-disc} in terms of the interpolants as
\begin{align}\nonumber
&\forall\,\wt{\W}\!\in\!V_h
:\ \ \ \
\mathcal{R}_1^{}\big(\underline{\Z}_{\rmn,\tau h};\wt{\W}_\rmt{-}\underline{\Z}_{\rmt,\tau h}\big)
+\Big\langle\mathcal{E}_{_\text{C}}'
\Big(\frac{\TAU\barw_{\rmn,\tau h}{+}\chi\underline\Z_{\rmn,\tau h}}
{\TAU+\chi}\Big),\wt{\W}_\rmn{-}\barw_{\rmn,\tau h}\Big\rangle
\\[-.2em]&\nonumber\hspace*{6em}
+\big\langle\overline{\P}_{\tau h},\wt{\W}{-}\barw_{\tau h}\big\rangle
\ge\mathcal{R}_1^{}\big(\underline{\Z}_{\rmn,\tau h};\barw_{\rmt,\tau h}{-}\underline{\Z}_{\rmt,\tau h}\big)
+\big\langle[\barF_\TAU^{}]_\W',\wt{\W}-\barw_{\tau h}\big\rangle
\\[.2em]&\hspace*{8em}\text{ with the traction }
\overline{\P}_{\tau h}=
\calP_h(\overline{\wt{\DIR}}_{\Dir,\tau},\barf_\tau,\barw_{\tau h})
\ \text{ for all }t\in[0,T].
\label{Biot-disc-g-h-disc+}
\end{align}
Replacing $\wt{\W}{-}\underline{\Z}_{\tau h}$ with $\wt{\Z}$, we obtain
\begin{multline}\label{Eq_Biot-disc-g-h+-}
\mathcal{R}_1^{}\big(\underline{\Z}_{\rmn,\tau h};\wt{\Z}_\rmt\big)
+\Big\langle\mathcal{E}_{_\text{C}}'
\Big(\frac{\TAU\barw_{\rmn,\tau h}{+}\chi\underline\Z_{\rmn,\tau h}}
{\TAU+\chi}\Big),\wt{\Z}_\rmn{+}\underline{\Z}_{\rmn,\tau h}{-}\barw_{\rmn,\tau h}\Big\rangle
+\big\langle\overline{\P}_{\tau h},\wt{\Z}{+}\underline{\Z}_{\tau h}{-}\barw_{\tau h}\big\rangle\\[-.2em]
\ge\mathcal{R}_1^{}\big(\underline{\Z}_{\rmn,\tau h};\barw_{\rmt,\tau h}{-}\underline{\Z}_{\rmt,\tau h}\big)
+\big\langle[\barF_\TAU^{}]_\W',
\wt{\Z}{+}\underline{\Z}_{\tau h}-\barw_{\tau h}\big\rangle.
\end{multline}
To facilitate the convergence proof, we note that
$\frac{\TAU}{\TAU+\chi}\barw_{\tau h}{+}\frac{\chi}{\TAU+\chi}\underline\Z_{\tau h}=\barz_{\tau h}$, cf.~\eqref{Biot-disc-g-h++},
and $\barw_{\tau h}{-}\underline{\Z}_{\tau h}=(\TAU{+}\chi)\DT\Z_{\tau h}$,
and also $\barw_{\rmt,\tau h}=\barz_{\tau h}{+}\chi\DT\Z_{\tau h}$,
cf.~\eqref{Biot-disc-g-h++}. Substituting it into \eqref{Biot-disc-g-h-disc+}
and realizing the 1-homogeneity of $\mathcal{R}_1^{}$ so that the factor
$\TAU{+}\chi$ can be forgotten, we obtain the discrete problem formulated
in terms of $\Z$'s only:
\begin{subequations}\label{Biot-disc-g-h-disc++}\begin{align}\nonumber
&\!\!\forall\,\wt{\Z}\!\in\!V_h:\ \ \
\mathcal{R}_1^{}\big(\underline{\Z}_{\rmn,\tau h};\wt{\Z}_\rmt\big)
+\big\langle\mathcal{E}_{_\text{C}}'
\big(\barz_{\rmn,\tau h}\big),\wt{\Z}_\rmn{-}\DT{\Z}_{\rmn,\tau h}\big\rangle
+\big\langle\overline{\P}_{\tau h},\wt{\Z}{-}\DT{\Z}_{\tau h}\big\rangle
\\&\label{Biot-disc-g-h-disc+++}\hspace*{17em}
\ge\mathcal{R}_1^{}\big(\underline{\Z}_{\rmn,\tau h};\DT{\Z}_{\rmt,\tau h}\big)
+\big\langle[\barF_\TAU^{}]_\W',\wt{\Z}-\DT{\Z}_{\tau h}\big\rangle
\\[.2em]&\label{Biot-disc-g-h-disc++++}\hspace*{5em}\text{ with the traction }
\overline{\P}_{\tau h}=
\calP_h(\overline{\wt{\DIR}}_{\Dir,\tau},
\barf_\tau,\barz_{\tau h}{+}\chi\DT{\Z}_{\tau h})
\ \text{ for all }t\in[0,T],
\end{align}\end{subequations}
and where the physical dimension of the test function $\wt{\Z}$ is, in contrast to the previous occurrence,  the same as for velocity due to the omitted factor $\TAU{+}\chi$.
An important additional information we can read from \eqref{Biot-disc-g-h-disc+++} is that
\begin{align}
\overline{\P}_{\tau h}\in\partial_{\DT\Z}\mathcal{R}_1^{}
\big(\underline{\Z}_{\rmn,\tau h};\DT{\Z}_{\rmt,\tau h}\big)-\mathcal{E}_{_\text{C}}'
\big(\barz_{\rmn,\tau h}\big)+[\barF_\TAU^{}]_\W'.
\end{align}
As the normal part of  $\partial_{\DT\Z}\mathcal{R}_1^{}$ is simply 0, we have
also an additional a-priori estimate of $[\overline{\P}_{\tau h}]_\rmn^{}=
[[\barF_\TAU^{}]_\W'
-\mathcal{E}_{_\text{C}}'(\barz_{\rmn,\tau h})]_\rmn^{}$.
As for the tangential part, we use that $|\partial_{\DT\Z}\mathcal{R}_1^{}|$ is
bounded by $\mu|\mathcal{E}_{_\text{C}}'(\barz_{\rmn,\tau h})|$, so that we have
$[\overline{\P}_{\tau h}]_\rmt^{}$ bounded but only in $L^\infty(0,T;L^{q_1}(\GC))$ for
any $q_1<\infty$; note that we cannot expect any bound in
$H^1(0,T;H^{1/2}(\GC))$ because $[\partial_{\DT\Z}\mathcal{R}_1^{}]_\rmt^{}$ is
truly set-valued. Altogether, we have $\overline{\P}_{\tau h}$ bounded in
$L^\infty(0,T;L^{q_1}(\GC;\R^d))$ for all
$1\le q_1<\infty$ (for $d=2$) or $1\le q_1\le4(q{-}1)$ (for $d=3$) with
$q$ an exponent from the growth condition \eqref{e4:friction-small-gamma}.

We want to pass to the limit in (the integral version of) the inequality
\eqref{Biot-disc-g-h-disc++}.
We can, for a while, select a subsequence and $(\Z,\P)\in H^1(0,T;H^{1/2}(\GC;\R^d))\times L^\infty(0,T;L^q(\GC;\R^d))$ such that
\begin{subequations}\label{conv}\begin{align}\label{conv1}
&\Z_{\tau h}\weak\Z&&\text{ in }\
H^1(0,T;H^{1/2}(\GC;\R^d)),
\\&\barz_{\tau h}\weaks\Z\ \text{ and }\ \underline{\Z}_{\tau h}\weaks\Z\!\!\!\!\!\!\!\!\!\!&&\text{ in }\ L^\infty(0,T;H^{1/2}(\GC;\R^d)),
\\\label{conv-p}
&\overline{\P}_{\tau h}\weaks\P&&\text{ in }\
L^\infty(0,T;L^{q_1}(\GC;\R^d))
\intertext{where $q_1$ is as above and, as usual, `$\weak$' or `$\weaks$'
mean convergence in
the weak or the weak* topology, respectively. By the (generalized)
Aubin-Lions theorem, we have also a strong convergence}
\label{conv-of-traces}
&\barz_{\tau h}\to\Z\ \text{ and }\ \underline{\Z}_{\tau h}\to\Z\!\!\!\!\!\!\!\!\!\!&&\text{ in }\ L^{q_0}(0,T;L^{q_2}(\GC;\R^d)),
\end{align}\end{subequations}
for all $1\le q_0<\infty$ and all $1\le q_2<\infty$ (for $d=2$) or $1\le q_2<4$
(for $d=3$).

The approximation $\calP_h$ of the Poincar\'e-Steklov operator
$\calP$  involved in \eqref{Biot-disc-g-h-disc++}
is to be specified. Standard ways of this approximation are by BEM
or by FEM. For the FEM-approximation of Poincar\'e-Steklov operators, see,
e.g., \cite{LanSte07CFBE}.
For a comparison with the BEM approximations of Poincar\'e-Steklov operators, see, e.g., \cite{GwiSte93BEPC,Pechstein2013,SauSch11BEM,MaiSte05hpBEM}.
Technically, the FEM-approximation is more amenable to convergence analysis because it smears a lot of technicalities and imitate essentially the conventional FEM arguments.

By this way, we can write the term $\langle\overline{\P}_{\tau h},\wt{\Z}{-}\DT{\Z}_{\tau h}\rangle$ in \eqref{Biot-disc-g-h-disc+++} with $\overline{\P}_{\tau h}$ from
\eqref{Biot-disc-g-h-disc++++} as a bulk integral, namely
\begin{equation}\label{p-vs-u}
\big\langle\overline{\P}_{\tau h},\wt{\Z}{-}\DT{\Z}_{\tau h}\big\rangle
=\sum_{\eta={\rm A,B}}\int_{\Omega^\eta}\!\!\bbC e(\baru_{\tau h}^\eta
{+}\chi\DT u_{\tau h}^\eta){:}e(\wt u^\eta{-}\DT u_{\tau h}^\eta)\,\dd\Omega
\end{equation}
where $\baru_{\tau h}\in H^1(\Omega\AA\cup\Omega\BB;\R^d)$ is the unique
solution to a finite-element approximation of the TBVP \eqref{BVP-k} with
$({\Vector{\DIR}}_\Dir,f\N,\W)$ taken as
$(\overline{\Vector{\DIR}}_{\Dir,\tau}{+}\DT\chi{\Vector{\DIR}}_{\Dir,\tau},
\barf_\tau,\barw_\tau)$, and $\wt{u}^\eta$ is an appropriate FEM-approximation test function which provides the interface gap $\wt{\Z}$.
Thus, similarly to the physical dimension of $\wt{\Z}$, used in~\eqref{Biot-disc-g-h-disc+++}, the dimension of $\wt{u}^\eta$ is also velocity.
After integration over $[0,T]$ by using the estimate
\begin{equation}\label{Eq_EstimBulkE}
\int_0^T\!\!\int_{\Omega^\eta}\bbC e(\baru_{\tau h}^\eta){:}
e(\DT u_{\tau h}^\eta)\,\dd\Omega
\ge\frac12\int_{\Omega^\eta}\bbC e(u_{\tau h}^\eta(T))
{:}e(u_{\tau h}^\eta(T))-\bbC e(u_0){:}e(u_0)\,\dd\Omega
\end{equation}
and a similar estimate for the normal-compliance term, we can write
\eqref{Biot-disc-g-h-disc++} in terms of $u$ in the form
\begin{align}\nonumber
&\sum_{\eta={\rm A,B}}\int_0^T\!\!\!\int_{\Omega^\eta}\!\!
\Big(\bbC e(\baru_{\tau h}^\eta){:}e(\widetilde u)+
\chi\bbC e(\DT u_{\tau h}^\eta)
{:}e(\widetilde u-\DT u_{\tau h}^\eta)\Big)\,\dd\Omega\dd t
-\!\int_0^T\!\!\!\int_{\GC}\!
\!\mu\gamma'\big(\JUMP{\underline u_{\tau h}}{\rmn}\big)
\big|\JUMP{\wt u}{\rmt}\big|\,\dd\Gamma\dd t
\\[-.3em]&\nonumber\qquad\qquad\qquad
+\int_{\GC}\!\!\gamma(\JUMP{u_0}{\rmn})\,\dd\Gamma
+\frac12\!\sum_{\eta={\rm A,B}}\int_{\Omega^\eta}\!\!\bbC e(u_0){:}e(u_0)\,\dd\Omega
\\&\nonumber\qquad
\ge\int_{\GC}\!\!\gamma\big(\JUMP{u_{\tau h}(T)}{\rmn}\big)\,\dd\Gamma
+\frac12\!\sum_{\eta={\rm A,B}}\int_{\Omega^\eta}\!\!\bbC e(u_{\tau h}^\eta(T))
{:}e(u_{\tau h}^\eta(T))\,\dd\Omega
\\[-.3em]&\qquad\qquad\qquad
-\int_0^T\!\!\!\int_{\GC}\!\!\mu\gamma'\big(\JUMP{\underline u_{\tau h}}{\rmn}\big)\big|\JUMP{\DT u_{\tau h}}{\rmt}\big|\,\dd\Gamma\dd t
+\sum_{\eta={\rm A,B}}\int_0^T\!\!\!\int_{\GN^\eta}f\N^\eta{\cdot}(\wt u-\DT u_{\tau h}^\eta)\,\dd\Gamma\dd t.
\label{conv2}
\end{align}
Then, we use
\begin{subequations}\label{conv-of-u}
\begin{align}\label{conv-of-u-strong}
&\baru_{\tau h}\to u&&\text{
 in } L^2(0,T;H^1(\Omega\AA{\cup}\Omega\BB;\R^d)),
\\&\DT u_{\tau h}\weak \DT u&&\text{
in }\ L^2(0,T;H^1(\Omega\AA{\cup}\Omega\BB;\R^d)),
\\&u_{\tau h}(T)\weak u(T)&&
\text{
 in }\ H^1(\Omega\AA\cup\Omega\BB;\R^d)),
\\\label{conv-un}&\JUMP{\underline u_{\tau h}}{\rmn}\to\JUMP{
u}{\rmn}
&&\text{
in }L^{q_0}(0,T;L^{q_1}(\GC;\R^{d-1})),\text{ and}
\\&\JUMP{\DT u_{\tau h}}{\rmt}\to\JUMP{\DT u}{\rmt}&&\text{
in }
L^2(0,T;L^{q_1}(\GC;\R^{d-1}));
\end{align}\end{subequations}
cf.~\eqref{conv-of-traces} for \eqref{conv-un}.
The strong convergence \eqref{conv-of-u-strong} is due to the uniform convexity
of the stored energy, counting that the nonmonotone friction term is in a
position of the lower-order term;  e.g.~\cite[Sect.\,5.2.4]{MieRou15RIS}
for details.
Thus, we can pass to the limit in the left hand side of \eqref{conv2} and
estimate from above the limit superior of the bulk terms there, while the
limit inferior of the right-hand side of \eqref{conv2} can be estimated
from below. Altogether, we thus obtain
\begin{align}\nonumber
&\sum_{\eta={\rm A,B}}\int_0^T\!\!\!\int_{\Omega^\eta}\!\!
\Big(\bbC e(u^\eta){:}e(\widetilde u)
+\chi\bbC e(\DT u^\eta)
{:}e(\widetilde u{-}\DT u^\eta)\Big)\,\dd\Omega
-\int_0^T\!\!\!\int_{\GC}\!
\!\mu\gamma'\big(\JUMP{u}{\rmn}\big)
\big|\JUMP{\wt u}{\rmt}\big|
\,\dd\Gamma
\\[-.3em]&\nonumber\qquad\qquad\qquad
+\int_{\GC}\!\!\gamma\big(\JUMP{u_0}{\rmn}\big)\,\dd\Gamma
+\frac12\!\sum_{\eta={\rm A,B}}\int_{\Omega^\eta}\!\!\bbC e(u_0){:}e(u_0)\,\dd\Omega
\\&\nonumber\qquad
\ge\int_{\GC}\!\!\gamma\big(\JUMP{u(T)}{\rmn}\big)\,\dd\Gamma
+\frac12\!\sum_{\eta={\rm A,B}}\int_{\Omega^\eta}\!\!\bbC e(u^\eta(T))
{:}e(u^\eta(T))\,\dd\Omega
\\[-.3em]&\qquad\qquad\qquad
-\int_0^T\!\!\!\int_{\GC}\!\!\mu\gamma'\big(\JUMP{u}{\rmn}\big)
\big|\JUMP{\DT u}{\rmt}\big|\,\dd\Gamma
+\sum_{\eta={\rm A,B}}\int_0^T\!\!\!\int_{\GN^\eta}f\N^\eta{\cdot}(\wt u{-}\DT u^\eta)\,\dd\Gamma
\label{conv2+}
\end{align}
for any $\wt u\in L^2(0,T;H^1(\Omega\AA\cup\Omega\BB))$. Actually,
\eqref{conv2+} is a weak formulation of the initial-boundary-value
problem \eqref{BVP}. As $u(t){+}\chi\DT u(t)$ solves the TBVP \eqref{BVP-k},
we can define the traction $p(t)$, and substitute the bulk integrals in
\eqref{conv2+} analogously like we did in the discrete form in \eqref{p-vs-u}.
In this way, we identify the limit traction $p$ from \eqref{conv-p} and we
arrive to
a weak formulation defined exclusively on $\GC$ in the form
\eqref{Biot-disc-g-h-disc++} written in the limit by omitting $\tau$ and $h$.
Eventually, putting $\wt u=0$ into \eqref{conv2+}, we obtain the upper energy
estimate
\begin{align}\nonumber
&\frac12\!\sum_{\eta={\rm A,B}}\int_{\Omega^\eta}\!\!\bbC e(u^\eta(T))
{:}e(u^\eta(T))\,\dd\Omega
+\int_{\GC}\!\!\gamma\big(\JUMP{u(T)}{\rmn}\big)\,\dd\Gamma
\\[-.5em]&\nonumber\qquad\qquad
+\chi\!\!\sum_{\eta={\rm A,B}}\int_0^T\!\!\!\int_{\Omega^\eta}\!\!\bbC e(\DT u^\eta)
{:}e(\DT u^\eta)\,\dd\Omega
-\int_0^T\!\!\!\int_{\GC}\!\!\mu\gamma'\big(\JUMP{u}{\rmn}\big)
\big|\JUMP{\DT u}{\rmt}\big|\,\dd\Gamma
\\&
\qquad
\le\int_{\GC}\!\!\gamma\big(\JUMP{u_0}{\rmn}\big)\,\dd\Gamma
+\frac12\!\sum_{\eta={\rm A,B}}\int_{\Omega^\eta}\!\!\bbC e(u_0){:}e(u_0)\,\dd\Omega
+\sum_{\eta={\rm A,B}}\int_0^T\!\!\!\int_{\GN^\eta}f\N^\eta{\cdot}\DT u^\eta\,\dd\Gamma.
\label{upper-engr-est}
\end{align}

\begin{remark}[{\sl Energy conservation.}]\label{rem-energy}
\upshape
As we have here linear visco-elastic material, one can modify the arguments
from \cite[Remark A.2]{KrPaRo14QACD}
to show even the equality in \eqref{upper-engr-est}.
\end{remark}

\begin{remark}[{\sl Vanishing viscosity towards rate-independent limit.}]
\label{rem-VV}
\upshape
Frequently, the relaxation time $\chi>0$ is very small with respect to the
external-loading time-scale, and the material is nearly elastic. However,
even very slow loading can lead to fast jumps during which the even very small
viscosity plays an essential role, as
demonstrated e.g. in~\cite{roubicek13A4} in case of an adhesive but frictionless contact.
Here, although in specific cases the numerical experiments may indicate
certain convergence as in Fig.\,\ref{Fig_SkewEFTR}, theoretical
justification hardly can be expected because the results for $\chi=0$
(i.e.\ quasistatic rate-independent friction problem) are missing
especially if the friction coefficient is not (unspecifically) small, in
spite of a long-lasting analytical effort in this direction worldwide.
\end{remark}

\end{document}